\documentclass[11pt]{article}
\usepackage{mathrsfs,amsfonts,amsmath}
\RequirePackage{ifthen,calc}
\usepackage{theorem}
\usepackage[dvips]{color}
\usepackage{graphicx}

\def\R{{\mathbb R}}

\def\Re {{\rm Re}\,}
\def\E{{\mathbb E}}
\def\P{{\mathbb P}}

\def\C{{\mathbb C}}

\newcommand{\red}[1]{\textcolor{red}{#1}}

\catcode`@=12

\newtheorem{thm}{\noindent Theorem}[section]
\newtheorem{lem}{\noindent Lemma}[section]
\newtheorem{cor}{\noindent Corollary}[section]

\newtheorem{remark}{Remark}[section]
\newtheorem{definition}{Definition}[section]
\theoremheaderfont{\normalfont\bfseries}
\theorembodyfont{\slshape} {\theorembodyfont{\rmfamily}} {\theorembodyfont{\rmfamily}
}
{\theorembodyfont{\rmfamily} } {\theorembodyfont{\rmfamily}
}
{\theorembodyfont{\rmfamily} } {\theorembodyfont{\rmfamily}
}
{\theorembodyfont{\rmfamily} } {\theorembodyfont{\rmfamily} } {\theorembodyfont{\rmfamily}
}

 \def\beqlb{\begin{eqnarray}}\def\eeqlb{\end{eqnarray}}
 \def\beqnn{\begin{eqnarray*}}\def\eeqnn{\end{eqnarray*}}
 \def\Re{\textrm{\it Re}}
 \numberwithin{equation}{section}
\def\qed{\hfill$\square$\smallskip}

 \newcommand{\pend}{\hfill \thicklines \framebox(5.5,5.5)[l]{}}
 \newenvironment{proof}{\noindent {\sc  Proof.} \rm}{\pend}

\begin{document}
\title{Stationary Distributions for Two-Dimensional Sticky Brownian Motions:  Exact Tail Asymptotics and Extreme Value Distributions  }
\author{Hongshuai Dai$^{a,b}$ and Yiqiang Q. Zhao$^b$
 \\ {\small a. School of Statistics, Shandong University of Finance and Economics, Jinan, Shandong
250014, P.R. China}\\
{\small b. School of Mathematics and Statistics,
Carleton University, Ottawa, ON, Canada K1S 5B6}}
\maketitle
\begin{abstract}
In this paper, we consider a two-dimensional sticky Brownian motion. Sticky Brownian motions can be viewed  as time-changed semimartingale reflecting Brownian motions, which find applications in many areas including queueing theory and mathematical finance.  For example,  a sticky Brownian motion can be used to model a storage system.with exceptional services. In this paper, we focus on stationary distributions for sticky Brownian motions. The main results obtained here include tail asymptotic properties in boundary stationary distributions, marginal distributions, and joint distributions.  The kernel method, copula concept and extreme value theory are main tools used in our analysis.
\end{abstract}

\small {{\bf MSC(2000):}  60K25, 60J10\\[1mm]


\section{Introduction}

Stochastic processes with sticky points in the Markov process sense have been studied extensively. A sticky Brownian motion on the half-line is the process evolving as a standard Brownian motion away from zero and reflecting at zero after spending a random time there --- as opposed to the one-dimensional semimartingale reflecting Brownian motion (SRBM), which reflects instantaneously. This process was initially studied by Feller \cite{F1952,F1954,F1957}, and It\^{o} and McKean \cite{IM1963, IM1965} in a more general context, and subsequently analyzed in more detail by many other authors \cite{B1999,IW1989}. These papers show that sticky Brownian motions arise as a time change of a reflecting Brownian motion, and that it describes the scaling limit of random walks on the
natural numbers whose jump rate at zero is significantly smaller than that at positive sites. Recently, R\'acz and Shrocnikov \cite{RS2015}  introduced multidimensional sticky Brownian motions which are an natural multidimensional  extension of sticky Brownian motions on the half-line.  As shown in \cite{RS2015}, a multidimensional sticky Brownian motion can also be written as a time-changed multidimensional SRBM.
Multidimensional sticky Brownian motions are of interest in both queuing theory and mathematical finance. For example, we can use it to model a market, which experiences a slowdown due to a major event (such as a court trial between some of the largest firms in the market) deciding about the new market leader, or a queueing system, in which the service time of the first customer (to an empty system) is different from other service times.
In the setting for single server queues, Welch~\cite{W1964} introduced an exceptional service for the first customer in each busy period and showed that a sticky Brownian motion on the half-line can be a heavy traffic limit. Later, with different exceptional service mechanisms, the same heavy traffic limit, or the skicy Brownian motion, was confirmed for other single server queueing models by Lemoine \cite{L1975},  Harrison and Lemoine \cite{HL1981}, Yamada \cite{Y1994}, and Yeo \cite{Y1961}. Similar to the one dimensional setting case, we expect applications in multi-server queueing networks.

 For a stable process, it is interesting and important to study its stationary probabilities. However, except for very limited special cases, we cannot get a closed-form solution for the stationary probability distribution. This adds values to study tail asymptotic properties in stationary probabilities, since performance bounds and approximations can often be developed from the tail asymptotic property. The stationarity of SRBMs has been studied in the literature. For example, Harrison and Hasenbein \cite{HH2009} and Harrison and Williams \cite{HW1987} studied the existence and uniqueness of the stationary distributions and  Franceschi and  Raschel~\cite{FR2017} obtained Laplace transforms for the defined boundary measures.  On the other hand, exact tail asymptotics for SRBM has been obtained recently, including Dai and Miyazawa~\cite{DM2011,DM2013}, Franceschi and Kurkova~\cite{FK2016}, and  Dai, Dawson and Zhao~\cite{DDZ}, in which tail asymptotic properties for the stationary distribution along a direction, or for the boundary measures were obtained.

In this paper, we extend the above research to study exact tail asymptotic properties for a time-changed SRBM. In additional to the asymptotic analysis for the boundary measures and for the distribution along a given direction, we also study exact tail asymptotic behaviours for the joint distributions. This research, was in part inspired by the recent work on time-changed SRBMs in R\'acz and Shrocnikov \cite{RS2015}. In \cite{RS2015}, the authors studies existence and uniqueness for the stationary distributions of multidimensional sticky Brownian motions. Furthermore, under a rather strict condition, they  \cite{RS2015} presented an explicit expression of the stationary distribution.

The main contributions made in this paper include, under a stability condition for a two-dimensional time-changed SRBM:
  \begin{itemize}
 \item[(1)] Exact tail asymptotics for the boundary stationary distributions, the marginal distributions, and the stationary distribution along a given direction; and
\item[(2)] The extreme value distribution and exact tail aymptotics for the joint stationary distribution.
  \end{itemize}

The rest of this paper is organized as follows:  In Section~\ref{sec:2}, we first recall some preliminaries related to sticky Brownian motions.  Section~\ref{sec:3} is devoted to studying basic properties of stationary distributions of the sticky Brownian motion. In Section~\ref{sec:4}, we apply the kernel method to study the tail behaviour for the boundary stationary distributions, the marginal stationary distributions and the joint distribution along a direction. In Section~\ref{sec:5}, we use copula concept and extreme value theory to study the tail behaviour of the joint stationary distribution. A final note is presented at the end of this paper.

\section{Preliminaries} \label{sec:2}

In this section, we introduce some preliminaries related to multidimensional sticky Brownian motions. We first recall the definition  of the semimartingale reflecting Brownian motion (SRBM). SRBM models arise as an approximation for queueing networks of various kinds (see for example, Williams~\cite{W1995,W}). A $d$-dimensional SRBM
$\tilde{Z}=\{\tilde{Z}(t),t\geq 0 \}$ is defined by the following:
 \beqlb\label{def1}
 \tilde{Z}(t)=X(t)+R L(t), \;\textrm{for}\;t\geq 0,
 \eeqlb
where  $\tilde{Z}(0)=X(0)\in\R^d_+$,  $X$ is an unconstrained Brownian motion with drift vector
$\mu=(\mu_1,\mu_2)'$ and covariance matrix $\Sigma=(\Sigma_{i,j})_{d\times d}$, $R=(r_{ij})_{d\times d}$ is a $d \times d$ matrix specifying the reflection behaviour at the boundaries,  and
$L=\{L(t)\}$ is a $d$-dimensional process with components $L_1, \ldots,
L_d$ such that:
\begin{itemize}
\item[(i)] $L$ is continuous and non-decreasing with $L(0)={\bf 0}$;

\item[(ii)] $L_j$ only increases at times $t$ for which
$\tilde{Z}_j(t)=0$, $j=1,\ldots, d$;

\item[(iii)] $\tilde{Z}(t)\in\R^d_+$, $t\geq 0$.
\end{itemize}

The existence of SRBM has been studied extensively, for example, Taylor and Williams~\cite{TW1993}, and
Reiman and Williams~\cite{RW1988}. Recall that  a $d\times d$ matrix $R$ is called an
$\mathbb{S}$-matrix, if there exists a $d$-vector $\omega\geq0$
such that $R\omega\geq0$, or equivalently, if there exists
$\omega>0$ such that $R\omega>0$. Furthermore, $R$ is called
completely $\mathbb{S}$ if each of its principal sub-matrices is
an $\mathbb{S}$-matrix. It was proved in \cite{TW1993, RW1988} that for a given set of
data $(\Sigma, \mu,R)$ with $\Sigma$ being positive definite,
there exists an SRBM for each initial distribution of $\tilde{Z}(0)$ if
and only if $R$ is completely $\mathbb{S}$. Furthermore, when $R$
is completely $\mathbb{S}$, the SRBM is unique in distribution for
each given initial distribution. It is well-known that a necessary condition (see, for example, Harrison and Williams~\cite{HW1987}, or Harrison and Hasenbein~\cite{HH2009}) for the
existence of the stationary distribution for $\tilde{Z}$ is
\beqlb\label{M-1} R \;\textrm{is non-singular and }\;R^{-1}\mu<0.
\eeqlb

\begin{remark}\label{cor-1}
When $d=2$, equation \eqref{M-1} is equivalent to
\beqlb\label{1-1}
    r_{11}>0,\; r_{22}>0, \; r_{11}r_{22}-r_{12}r_{21}>0; \eeqlb
and
\beqlb\label{1-2}
    r_{22}\mu_1-r_{12}\mu_2<0,\; r_{11}\mu_2-r_{21}\mu_1<0.
\eeqlb
\end{remark}

We note that SRBM does not spend time on the boundary. As opposed to it, a sticky Brownian motion would spend a duration of  time on the boundary.  For the one-dimensional case, Feller~\cite{F1952,F1954,F1957} first observed the sticky boundary behaviour for diffusion processes, and studied the problem that describes domains of the infinitesimal generators associated with a strong Markov process $\tilde{X}$ in $[0,\,\infty)$.  Moreover, $\tilde{X}$ behaves like a standard Brownian motion  in $(0,\;\infty)$, while at $0$, a possible boundary behaviour is described by
\beqlb\label{2-a1}
f'(0+)=\frac{1}{2u}f''(0+)
\eeqlb
where $u\in(0,\infty)$ is a given and fixed constant, and $f$ are functions belonging to the domain of the infinitesimal generator of $\tilde{X}$.
The second derivative $f''(0+)$ measures the ``stickiness" of $\tilde{X}$ at $0$. For this reason, the process $\tilde{X}$ is called a sticky Brownian motion (sometimes, it is also referred to as a sticky reflecting Brownian motion in the literature). It\^{o} and Mckean~\cite{IM1963} first constructed the sample paths of $\tilde{X}$.  They showed that $\tilde{X}$ can be obtained from one dimensional SRBM $\tilde{Z}$ by the time-change $t\to T(t):=S^{-1}(t)$,  where $S(s)=s+\frac{1}{u}L_s$ for $s>0$, or $T(t)=s$ is determined by the equation $t = s+\frac{1}{u}L_s$.
For more information about sticky Brownian motions on the half-line, refer to Engelbert and Peskir~\cite{EP2014} and references therein.

 R\'acz and Shrocnikov~\cite{RS2015} introduced multidimensional sticky Brownian motions, which is a natural extension of the sticky Brownian motion on the half-line, and proved existence and uniqueness of the multidimensional sticky Brownian motion.  Similar to a sticky Brownian motion on the half-line,  let
\beqlb\label{2-a3}
S(t)=t+\sum_{i=1}^du_iL_i(t),
\eeqlb
where $u_i\in(0,\;\infty)$, $i=1,\ldots,d$, and let $T(\cdot)$ be the inverse of $S(t)$, that is,
 $T(t)=S^{-1}(t)$ .
\begin{remark}\label{rem-1}
It  follows from Kobayashi \cite{K2011} and \eqref{2-a3} that  $T$ has continuous paths and $\lim_{t\to\infty}T(t)=\infty.$  Moreover,  $T(1)\leq 1$.
  \end{remark}

According to the above discussion, a multidimensional sticky Brownian motion can be defined as
\beqlb\label{1-4}
Z(t)=\tilde{Z}(T(t)).
\eeqlb
This type of processes finds applications in the fields of queueing theory and finance. In the queueing field,  it is well known that SRBM is a heavy traffic limit for many queuening networks such as open queueing networks. As discussed in the introduction, in the setting for single server queues, a sticky Brownian motion on the half-line can be served as a heavy traffic limit of a queueing system with exceptional service mechanisms. It is reasonable to expect that a multidimensional sticky Brownian motion serves as a heavy traffic limit for such multidimensional queueing networks with appropriately defined exceptional service mechanisms.

\section{Basic Adjoint Relation} \label{sec:3}

Establishment of the basic adjoint relation (BAR) is the starting point for the analysis of our work in this paper, which is the counterpart to the fundamental form in the discrete case. Based on this equation, we can extend the kernel method (for example, see Li and Zhao~\cite{LZ2015} and references therein) to study exact tail asymptotics for stationary distributions of a sticky SRBM.  This is the focus of this section.

In the rest of this paper, we assume that $Z(0)$ follows the stationary distribution $\pi$ of $Z(t)$. Recall that we can define the moment generating function (MGF) of any finite measure on $\mathscr{B}(\R_+^2)$, where $\R_+^2= \{ x=(x_1,x_2): x_i \geq 0 \}$. For example, for the stationary measure $\pi$, the MGF $\Phi(\theta)$ is defined as follows:
\beqnn
    \Phi(\theta)=\int_{\R_+^2}\exp\{<\theta,x>\} \pi(dx).
\eeqnn
We study tail asymptotic properties in $\pi$ through the kernel method by analyzing the kernel function in the BAR, which connects $\Phi(\theta)$ to other unknown MGFs, for which analysis is possible.

Similar to SRBM, $\Phi(\theta)$ is closely related to the MGFs of two boundary measures, which are defined below.
For any set $A\in\mathscr{B}(\R_+^2)$, define
\beqlb\label{2-25}
V_i(A)&&=\E_{\pi}\bigg[\int_0^{T(1)} 1_{\{\tilde{Z}(s)\in A\}}dL_i(s)\bigg].
\eeqlb
In addition, for any Boreal measure $B\in\mathscr{B}(\R_+^2)$, we define the joint measure for the time change:
\beqnn
    V_0(B)=\E_{\pi}\bigg[\int_0^1 1_{\{Z(s)\in B\}}dT(s)\bigg].
\eeqnn
According to Corollary~\ref{cor-2} below, all $V_i$, $i=0, 1,2$, are finite measures on $\R_+^2$.
Then, we can define MGFs $\Phi_i(\theta)$ for $V_i$, $i=0,1,2,$ by
\beqnn
\Phi_i(\theta)=\int_{\R_+^2}\exp\{<\theta,x>\}V_i(dx).
\eeqnn

For these measures,  we have the following BAR.
\begin{lem}\label{thm-1}\quad
\begin{itemize}
\item[(1)]The boundary measures $V_i$, $i=1,2,$ and the joint measure $V_0$ are all finite.
\item[(2)]The MGFs of $V_i$, $i=0,1,2,$ have the following BAR: for any ${\bf\theta}\in \R_{-}^2 = \{ \theta = (\theta_1, \theta_2)': \theta_i < 0\}$,
\beqlb\label{2-5-a}
-\Psi_{X}({\bf\theta})\Phi_0({\bf\theta})=\Phi_1({\bf\theta})<{\bf\theta}, R_1>+\Phi_2({\bf\theta})<{\bf\theta}, R_2>,
\eeqlb
where $R_i$ is the $i$th column of the reflection matrix $R$, and $\Psi_{X}(\theta)$ is the L\'evy exponent of the multidimensional Brownian vector $X(1)$.
\end{itemize}
\end{lem}

\proof Since $Z(0)$ follows the stationary distribution $\pi$, for any $t\in\R_+$,
\beqnn\label{2-27}
\P(Z(t)\leq z)=\P(Z\leq z).
\eeqnn

We note that $\{Z(t)\}$ is a semimartingale. Since $T(t)$ is continuous and $S(t)$ is strictly increasing,  it follows from Kobayashi~\cite[Corollary 3.4]{K2011} that  if $f:\R^2\to \R$ is $C_{b}^2$ function, then
\beqlb\label{3-5}
&&f(Z(t))-f(Z(0))=\sum_{i=1}^2 \mu_i\int_0^{T(t)}\frac{\partial f}{\partial x_i}(\tilde{Z}(u))du+ \sum_{i,j=1}^2 \int_0^{T(t)}r_{ji}\frac{\partial f}{\partial x_j}\big(\tilde{Z}(u)\big)dL_i(u)\nonumber
\\&&\hspace{1cm}+\sum_{i=1}^2\int_0^{T(t)}\frac{\partial f}{\partial x_i}\big(\tilde{Z}(u)\big)dX_i(u)+\frac{1}{2}\sum_{i,j=1}^2\Sigma_{i,j}\int_0^{T(t)}\frac{\partial^2 f}{\partial x_i\partial x_j}\big(\tilde{Z}(u)\big)du.
\eeqlb
Hence, we have
\beqlb\label{2-6}
&&\sum_{i=1}^2 \mu_i\E_{\pi}\bigg[\int_0^{T(1)}\frac{\partial f}{\partial x_i}(\tilde{Z}(u))du\bigg]+ \sum_{i,j=1}^2 \E_\pi\bigg[\int_0^{T(t)}r_{ji}\frac{\partial f}{\partial x_j}\big(\tilde{Z}(u)\big)dL_i(u)\bigg]\nonumber
\\&&\hspace{1cm}+\frac{1}{2}\sum_{i,j=1}^2\Sigma_{i,j}\E_\pi\bigg[\int_0^{T(t)}\frac{\partial^2 f}{\partial x_i\partial x_j}\big(\tilde{Z}(u)\big)du\bigg]=0.
\eeqlb
Next, we prove the firs part of this theorem.   From \eqref{2-25}, we get that for all $i=1,2$,
\beqlb
V_i(\R_+^2)=\E_\pi[ L_i\big(T(1)\big)],
\eeqlb
and
\beqnn
V_0(\R_+^2)=\E_\pi(T(1)).
\eeqnn
Hence, it suffices to prove that for any $i=1,2$,
\beqlb\label{2-28}
\E_\pi[ L_i\big(T(1)\big)]<\infty,
\eeqlb
and
\beqlb\label{2-28a}
\E_\pi(T(1))<\infty.
\eeqlb
Let $f(x,y)=\exp\{\theta x\}$ with $\theta<0$ and $x\geq 0$. Then we have that
\beqlb\label{2-32}
f'_1(x,y)=\theta\exp\{\theta x\}\;\text{and}\;f''_{1,1}(x,y)=\theta^2\exp\{\theta x\}.
\eeqlb
Hence, combing  \eqref{2-6} and \eqref{2-32} gives
\beqlb\label{2-31}
-\Psi_{X}(\theta,0)\Phi_0(\theta,0)=\E_\pi[L_1(T(1))]r_{11}\theta+ \Phi_2(\theta,0)r_{12}\theta.
\eeqlb
Dividing $\theta<0$ at both sides of \eqref{2-31} and letting $\theta\to 0$, we get that
\beqlb\label{2-33}
-\mu_1\E_{\pi}[T(1)]=\E_\pi[L_1(T(1))]r_{11}+\E_\pi[L_2(T(1))]r_{12}.
\eeqlb
Symmetrically, let $f(x,y)=\exp\{\theta y\}$ with $\theta<0$ and $y\geq 0$. Similar to \eqref{2-33}, we can get that
\beqlb\label{2-34}
- \mu_2\E_{\pi}[T(1)]=\E_\pi[L_1(T(1))]r_{21}+\E_\pi[L_2((T(1))]r_{22}.
\eeqlb
On the other hand, from the relation between $T(\cdot)$ and $S(\cdot)$, we get that
\beqlb\label{2-a2}
T(t)=t-\sum_{i=1}^2u_iL_i(T(t)).
\eeqlb
Hence, from \eqref{2-a2}, we have
\beqlb\label{2-aa2}
\E_{\pi}[T(1)]=1-\sum_{i=1}^2\E_{\pi}u_iL_i(T(1))
\eeqlb
Combing \eqref{2-33}, \eqref{2-34}  and \eqref{2-aa2} leads to \eqref{2-28} and \eqref{2-28a}.

Taking $f(x,y)=\exp\{\theta_1 x+\theta_2 y\}$ with $\theta_i\leq 0$, $i=1,2$, in equation \eqref{2-6} can prove the second part of this theorem. \qed

\begin{remark}  Let $ C^{2}_{b}(\R^{2}_{+})$  be the set of functions $f$ on $\R_+^2$  such that $f$, its first order derivatives, and its second order derivatives are
bounded and continuous. For any $f\in C^{2}_{b}(\R^{2}_{+})$, it follows from \eqref{2-6} that
\beqnn
\int_{\R_+^2}\mathcal{L}f(x) V_0(dx)+\sum_{i=1}^2\int_{\R_+^2}<\bigtriangledown f(x), R_{i}>V_i(dx)=0,
\eeqnn
where
\beqnn
\mathcal{L}f(x)=\frac{1}{2}\sum_{i,j=1}^2\Sigma_{i,j}\frac{\partial^2f}{\partial x_i\partial x_j}(x)+\sum_{j=1}^2\mu_j\frac{\partial f}{\partial x_j}(x),
\eeqnn
and $\bigtriangledown f(x)$ is the gradient of $f$.
From Dai and Kurtz~\cite[Theorem 1.4]{DK2003} (or Braverman, Dai and Miyazawa \cite[Lemma 2.1] {BDM2017}), we can get that $V_0(\cdot)/\E_\pi\big(T(1)\big)$, and $V_i(\cdot)/\E_\pi\big(L_i(T(1))\big)$, $i=1,2$, are the stationary distribution, and the boundary distribution of the corresponding reflecting Brownian motion $\tilde{Z}$, respectively.
\end{remark}
The following corollary immediately follows from the proof of Lemma~\ref{thm-1}.
\begin{cor}\label{cor-2}
\beqnn
\E\big[L_1(T(1))\big]=\frac{\mu_1(r_{22}-\mu_2u_2)-\mu_2(r_{12}-\mu_1u_2)}{(r_{21}-\mu_2u_1)(r_{12}-u_2\mu_1)-(r_{11}-\mu_1u_1)(r_{22}-\mu_2u_2)}
\eeqnn
and
\beqnn
\E\big[L_2(T(1))\big]=\frac{\mu_1(r_{21}-\mu_2u_1)-\mu_2(r_{11}-\mu_1u_1)}{(r_{22}-\mu_2u_2)(r_{11}-\mu_1u_1)-(r_{12}-\mu_2u_1)(r_{21}-\mu_2u_1)}.
\eeqnn
\end{cor}

Below, we state the main result of this section.  The sticky Brownian motion defined by \eqref{1-4} satisfies the following BAR.
\begin{thm}
\beqlb\label{2-5}
-\Psi_{X}({\bf\theta})\Phi({\bf\theta})=\Phi_1({\bf\theta})\big(<{\bf\theta}, R_1>-u_1\Psi_{X}({\bf\theta})\big)+\Phi_2({\bf\theta})\bigg(<{\bf\theta}, R_2>-u_2\Psi_{X}({\bf\theta})\bigg).
\eeqlb
\end{thm}

\proof
For any Boreal set $B\in\mathscr{B}(\R_+^2)$, we have
\beqlb\label{2-8}
\pi(B)&=&\E_\pi\Big[\int_0^11_{\{Z(s)\in B\}}ds\Big]\nonumber
\\&=&\E_\pi\Big[\int_0^{T(1)} 1_{\{\tilde{Z}(s)\in B \}}dS(s)\Big]\nonumber
\\&=& \E_\pi\Big[\int_0^{T(1)} 1_{\{\tilde{Z}(s)\in B \}}dt\Big]+\sum_{i=1}^2u_i\E_\pi\Big[\int_0^{T(1)} 1_{\{\tilde{Z}(s)\in B\}}dL_i(s)\Big]\nonumber
\\&=& V_0(B)+\sum_{i=1}^2 u_iV_i(B).
\eeqlb
From \eqref{2-5-a} and \eqref{2-8}, we can get \eqref{2-5}.\qed

\section{Exact Tail Asymptotics} \label{sec:4}

It is well-known that, except for some special cases, it is usually not expected to have an explicit expression for the stationary distribution for a multi-dimensional stochastic network. Instead, explicit tail asymptotic properties in the distribution could provide insightful understanding of the mode, and lead to performance bounds and numerical algorithms. Our focus in this section is on the so-called exact tail asymptotic behaviour for a tail probability function $g(x)$, which means to identify an explicitly expressed function $h(x)$ such that
\beqnn
\lim_{x\to\infty}\frac{g(x)}{h(x)}=1,
\eeqnn
denoted by $g(x) \sim h(x)$. A Tauberian-like theorem is used in the kernel method to link the asymptotic property for the tail probability function to the asymptotic property of the transform function of the distribution (see details, for example, in \cite{LZ2015}). To this end, we first study the kernel equation
\beqlb\label{3-a1}
    \Psi_{X}(x,y)=0.
\eeqlb

For $(x,y)$ satisfying \eqref{3-a1}, if the MGF $\Phi(x;y) < \infty$, then from  \eqref{2-5-a}, we have
\beqlb\label{3-2}
\gamma_1(x,y)\Phi_1(x,y)+\gamma_2(x,y)\Phi_2(x,y)=0,
\eeqlb
where $\gamma_1(x,y)=xr_{11}+y r_{21}$ and $\gamma_2(x,y)=xr_{12}+yr_{22}$.  Equation \eqref{3-2}
 provides a relationship between the two unknown MGFs $\Phi_i(x,y),i=1,2$. By using \eqref{3-2}, we can perform a singularity analysis of these functions to obtain exact tail asymptotics for the two boundary distributions $V_i, i = 1,2$, through the following steps:
\begin{itemize}
\item[(i)] Analytic continuation of the functions $\Phi_i(x,y)$ for $i=1,2$;
\item[(ii)] Singularity analysis of the functions $\Phi_i(x,y)$ for $i=1,2$; and
\item[(iii)] Applications of a Tauberian-like theorem, Theorem~\ref{TLT} below.
\end{itemize}

We can then further to obtain exact tail asymptotics for the marginal distributions and for the joint stationary distribution along a given direction.

To consider the analytic continuation of MGFs, we need the following version of Pringsheim's theorem for MGFs (see, for example, Dai and Miyazawa~\cite{DM2011} or Markushevich~\cite{M1977}).
\begin{lem}\label{2-lem-1}
Let $g(x)=\int_{0}^\infty e^{\lambda x}dF(x)$ be the moment generating function of a probability distribution $F$ on $\R_+$ with real variable $\lambda$. Define the convergence parameter of $g$ as
\beqnn
C_p(g)=\sup\{\lambda\geq 0: g(\lambda)<\infty\}.
\eeqnn
Then, the complex variable function $g(z)$ is analytic on $\{z\in \mathbb{C}: \Re (z)<C_p(g)\}$.
\end{lem}

For the bridge of the tail asymptotic behaviour in a stationary distribution and the asymptotic behaviour of its transform, we use the following Tauberian-like theorem, which can be found in Dai, Dawson and Zhao~\cite{DDZ}.
\begin{thm}\label{TLT} Let $g(s)$ be the Laplace-transformation of $f(s)$, i.e,
 \beqnn g(s)=\int_{0}^\infty e^{st}f(t)dt.
 \eeqnn
 Moreover, denote
 \beqnn
    \Delta(z_0,\epsilon)=\big\{z\in\C:z\neq z_0,\;|{\rm
arg}(z-z_0)|>\epsilon\big\},
 \eeqnn
 where ${\rm arg}(z) \in(-\pi,\;\pi]$ is the principal part of the argument of a complex number $z$.
Assume that $g(z)$ satisfies the following conditions:
\begin{itemize}
\item[(1)]The left-most singularity of $g(z)$ is $\alpha_0$ with
$\alpha_0>0$. Furthermore, we assume that as $z\to \alpha_0$,
$$g(z)\sim(\alpha_0-z)^{-\lambda}$$ for some
$\lambda\in\C\setminus\mathbb{Z}_{\leq 0}$;

\item[(2)]$g(z)$ is
analytic on $\mathbb{Z}_{< \alpha_0} \cap \Delta(\alpha_0,\epsilon_0)$ for some
$\epsilon_0>0$;

\item[(3)] $g(z)$ is bounded on
$\mathbb{Z}_{< \alpha_0} \cap \Delta(\alpha_0,\epsilon_1)$ for some $\epsilon_1>0$.
\end{itemize}
Then, as $t\to\infty$,
\beqnn
f(t)\sim e^{-\alpha_0 t}\frac{t^{\lambda-1}}{\Gamma(\lambda)},
\eeqnn
where $\Gamma(\cdot)$ is the Gamma function.
\end{thm}

In the rest of this paper, let $K$ denote an unspecified constant, whose value might be different from one case to another.

\subsection{Tail Asymptotics for Boundary Distributions}

In this subsection, we study exact tail asymptotics for the boundary probabilities $V_i$, $i=1,2$. This can be done by studying the kernel equation $\Psi_X(x,y)=0$ and the functions $\gamma_i(x,y)$, $i=1,2$, which is in parallel to the work in Dai, Dawson and Zhao~\cite{DDZ}. In \cite{DDZ}, the authors obtained exact tail asymptotics for the boundary stationary distributions associated to SRBM, while in this subsection the boundary distributions are associated with $Z$, a time-chaged SRBM. Since the similarity to \cite{DDZ} in the analysis, we provide the asymptotic results without proofs here.

 We first rewrite the kernel equation \eqref{3-a1} in a quadratic form in $y$
with coefficients that are polynomials in $x$:
 \beqlb\label{3-1}
\gamma(x,y):=\Psi_{X}(x,y)&=&x\mu_1+y\mu_2+\frac{1}{2}\Sigma_{11}x^2+\Sigma_{12}xy+\frac{1}{2}\Sigma_{22}y^2\nonumber
\\&=&\frac{1}{2}\Sigma_{22}y^2+(\mu_2+\Sigma_{12}x)y+\frac{1}{2}\Sigma_{11}x^2+x\mu_1\nonumber
\\&=& a y^2+b(x)y+c(x)=0,
\eeqlb where \beqnn
a=\frac{1}{2}\Sigma_{22},\;b(x)=\mu_2+\Sigma_{12}x
 \quad \textrm{and}\quad c(x)=x\mu_1+\frac{1}{2}\Sigma_{11}x^2.
 \eeqnn

Let
 \beqlb\label{4-13} D_1(x)=4\Big[(\Sigma_{12}^2-\Sigma_{11}\Sigma_{22})x^2+2(\Sigma_{12}\mu_1-\Sigma_{22}c_2)x+\mu_2^2\Big].
 \eeqlb
be the discriminant of the quadratic form in (\ref{3-1}).
Therefore, in the complex plane $\mathbb{C}_x$,  for every $x$,  two solutions to (\ref{3-1}) are given
by \beqlb\label{A-9} Y_{\pm}(x)=\frac{-b(x)\pm\sqrt{b^2(x)-4ac(x)}}{2a},
\eeqlb unless $D_1(x)=0$, for which $x$ is called a branch point
of $Y$.   We emphasize that in using the kernel method, all functions and variables are usually treated as
complex ones.

Symmetrically, when $x$ and $y$ are interchanged, we have
\beqlb\label{M-3}
\gamma(x,y)=\tilde{a}x^2+\tilde{b}(y)x+\tilde{c}(y)=0,
 \eeqlb
where
 \beqnn \tilde{a}=\frac{1}{2}\Sigma_{11},
\tilde{b}(y)=\Sigma_{12}y+\mu_1,\quad \textrm{and}\quad
\tilde{c}(y)=\frac{1}{2}\Sigma_{22}y^2+y\mu_2.
 \eeqnn

Let $D_2(y)=\tilde{b}^2(y)-4\tilde{a}\tilde{c}(y)$.
For each fixed $y$,  two solutions to (\ref{M-3}) are given by
\beqlb\label{M-4}
X_{\pm}(y)=\frac{-\tilde{b}(y)\pm\sqrt{\tilde{b}^2(y)-4\tilde{a}\tilde{c}(y)}}{2\tilde{a}},
\eeqlb
 unless $D_2(y)=0$, for which $y$ is called a branch point of $X$.

 By some basic calculations, we have the following result.
\begin{lem}\label{3-lem1}
$D_1(x)$  has two zeros satisfying $x_1\leq0<x_2$   with $x_i,i=1,2$, being real numbers.  Furthermore,
$D_1(x)>0$ in $(x_1,x_2)$, and $D_1(x)<0$ in $(-\infty,
x_1)\cup(x_2,\infty)$.  Similarly, $D_2(y)$ has two zeros
satisfying $y_1\leq0<y_2$   with $y_i,i=1,2$, being real numbers. Moreover, $D_2(y)>0$ in $(y_1,y_2)$,
and $D_2(y)<\infty$ in $(-\infty,y_1)\cup(y_2,\infty)$.
\end{lem}
 \begin{remark} It is obvious that
  \beqlb\label{4-4}
x_2=\frac{2(\Sigma_{12}\mu_2-\Sigma_{22}\mu_1)+\sqrt{\Delta}}{2(\Sigma_{11}\Sigma_{22}-\Sigma^2_{12})},
\eeqlb
where $\Delta=\big(\mu_2\Sigma_{12}-\Sigma_{22}\mu_1\big)^2-(\Sigma^2_{12}-\Sigma_{11}\Sigma_{12})\mu_2^2$.
 \end{remark}

 Next, we carry out the analytic continuations of MGFs $\Phi_i$, $i=1,2.$
 For convenience, denote
 \beqnn
    \tilde{\mathbb{C}}_x &=& \mathbb{C}_x\setminus\Big\{(-\infty, x_1]\cup
[x_2,\infty)\Big\}.
\\ \tilde{\mathbb{C}}_y &=& \mathbb{C}_y  \setminus\Big\{(-\infty, y_1]\cup[y_2,\infty)\Big\}.
 \eeqnn
Furthermore,  in the following we use $Y_0$ and $Y_1$ instead of
$Y_{-}$ and $Y_+$. Similarly, we use $X_0$ and $X_1$
instead of $X_{-}$ and $X_+$, respectively.
The following results in this subsection are obtained based on analysis of the BAR in (\ref{2-5-a}).

\begin{lem}\label{4-lem-5}
$\Phi_2(z):=\Phi_2(z,0)$ can be analytically continued to the region: $\{z \in \tilde{\mathbb{C}}_x: \gamma_2(z,Y_0(z))\neq 0\} \cap \{z \in \tilde{\mathbb{C}}_x: \Re (Y_0(z))<\tau_2\}$ with $\tau_2=C_p(\Phi_2)$, and
\beqlb\label{4-50}
\Phi_2(z)=-\frac{\gamma_1\big(z,Y_0(z)\big)\Phi_1\big(0,Y_0(z)\big)}{\gamma_2\big(z,Y_0(z)\big)}.
\eeqlb

Similarly, $\Phi_1(z):=\Phi_1(0,z)$,  can be analytically continued to the region: $\{z \in \tilde{\mathbb{C}}_y: \gamma_1(X_0(z),z)\neq 0\} \cap \{z \in \tilde{\mathbb{C}}_y: \Re (X_0(z))<\tau_1\}$ with $\tau_1=C_p(\Phi_1)$,
 and
\beqlb\label{M-5}
\Phi_1(z)=-\frac{\gamma_2\big(X_0(z),z\big)\Phi_2\big(X_0(z)\big)}{\gamma_1\big(X_0(z),z\big)}.
\eeqlb
\end{lem}

For a proof here, readers can refer to the proofs to Lemmas 3.2 and 4.4 in \cite{DDZ}.

\begin{lem}\label{4-lem7}
The function $\Phi_2(x)$ is  meromorphic on the cut plane $\tilde{\mathbb{C}}_x$. Similarly, the function $\Phi_1(y)$ is  meromorphic  on the cut plane $\tilde{\mathbb{C}}_y$.
\end{lem}

For a proof here, readers can refer to the proof to Lemma~10 in \cite{DDZ}.

Next, we study the tail behaviour of MGFs $\Phi_i$, $i=1,2$.  From Lemmas \ref{4-lem-5} and \ref{4-lem7}, we can obtain that the location of $\tau_2$ plays a vital role in determining the tail behaviours of $\Phi_2$.  Hence, we need to study the zeros of $\gamma_i(x,y)$, $i=1,2$, and compare it to the singularity $x_2$ to locate the dominant singularity of $\Phi_i$. For convenience, let $\tilde{x}$ be the solution of $\tilde{y}=Y_0(x)$ for $\gamma_1(X_0(\tilde{y}),\tilde{y})=0$, and $x^*$ be the zero of $\gamma_2\big(z,Y_0(z)\big)=0$. Based on the possibilities of the locations of $\tilde{x}$, $x^*$ and $x_2$,
we have the following results.
\begin{lem}\label{thm1}
For the function $\Phi_2(x)$, a total of four types of asymptotic properties exist as $x$ approaches to $\tau_2$, based on the detailed property of $\tau_2$.
\begin{description}
\item[Case 1:] If $\tau_2=x^*<\min\{\tilde{x},x_2\}$, or
$\tau_2=\tilde{x}<\min\{x^*,x_2\}$, or $\tau_2=\tilde{x}=x^*=x_2$,
then
    \beqlb\label{thm1-1}
    \lim_{x\to \tau_2}(\tau_2-x)\Phi_2(x)=A_1(\tau_2);
    \eeqlb

\item[Case 2:] If $\tau_2=x^*=x_2<\tilde{x}$, or
$\tau_2=\tilde{x}=x_2<x^*$, then \beqlb\label{thm1-2}
\lim_{x\to\tau_2}\sqrt{\tau_2-x}\Phi_2(x)=A_2(\tau_2); \eeqlb

\item[Case 3:] If $\tau_2=x_2<\min\{\tilde{x},x^*\}$, then
\beqlb\label{thm1-3}
\lim_{x\to\tau_2}\sqrt{\tau_2-x}\Phi_2'(x)=A_3(\tau_2); \eeqlb

\item[Case 4:] If $\tau_2=x^*=\tilde{x}<x_2$, then
\beqlb\label{thm1-4}
\lim_{x\to\tau_2}(\tau_2-x)^2\Phi_2(x)=A_4(\tau_2). \eeqlb
In the above, $A_i(\tau_2) \neq 0$, $i=1, 2, 3, 4$, is a constant depending on $\tau_2$.
\end{description}
\end{lem}

The proof of  Lemma~\ref{thm1} is based on Lemma~\ref{4-lem-5} and a  comparison between $x_2$, $x^*$ and $\tilde{x}$. For details, readers can refer to the proof to Theorem~1 in \cite{DDZ}.

Finally, we need to convert the asymptotic property, at the dominant singularity, of $\Phi_2$ to that for the tail probability of $V_2$ by the Tauberian-like theorem, Theorem~\ref{TLT}. One can check that all conditions in Theorem \ref{TLT} are satisfied by $\Phi_2(x)$.  Hence, we have

\begin{thm}\label{thm-a2}For the  boundary measure $V_2\big(\cdot)$, we have the following tail asymptotic properties for large $x$.
\begin{description}
\item[Case 1:] If $\tau_2=x^*<\min\{\tilde{x},x_2\}$, or
$\tau_2=\tilde{x}<\min\{x^*,x_2\}$, or $\tau_2=\tilde{x}=x^*=x_2$,
then
    \beqlb\label{4-9}V_2\big(x,\infty\big)\sim Ke^{-\tau_2 x};\eeqlb
\item[Case 2:] If $\tau_2=x^*=x_2<\tilde{x}$, or
$\tau_1=\tilde{x}=x_2<x^*$, then $$V_2\big(x,\infty\big)\sim Ke^{-\tau_2
x}x^{-\frac{1}{2}};$$

\item[Case 3:] If $\tau_2=x_2<\min\{\tilde{x},x^*\}$, then
$$V_2\big(x,\infty\big)\sim Ke^{-\tau_2
x}x^{-\frac{3}{2}};$$

\item[Case 4:] If $\tau_2=x^*=\tilde{x}<x_2$, then
$$V_2\big(x,\infty\big)\sim K e^{-\tau_2
x}x^{1}.$$
\end{description}

\end{thm}

In order to applying Theorem~\ref{thm-a2},  we need to determine the value of $\tau_2$. We have the following technical lemma, which can be useful.
\begin{lem}\label{5-lem2}\quad
\begin{itemize}
\item[(1)]$x^*$ is the root of $\gamma_2\big(z,Y_0(z)\big)=0$ in $(0,\;x_2]$ if and only if $\gamma_2\big(x_2,Y_0(x_2)\big)\geq 0$.
\item[(2)]If $\tilde{x}=\tau_2$, then $\tilde{x}=X_1(\tilde{y})$.
\end{itemize}
\end{lem}

\begin{remark}
 By some calculations, we can get

 \noindent
\beqlb\label{4-6}
\tilde{x}=
\left\{
\begin{array}{ll}
\displaystyle \frac{r_{21}}{r_{11}}y^*-2\frac{\Sigma_{12}y^*+\mu_1}{\Sigma_{11}},
&\;\textrm{if}\; Y_0(\tilde{x})= y^*; \\
\infty,&\;\textrm{otherwise},
\end{array}
\right.
\eeqlb
where
$$
y^*=\frac{2r_{11}r_{21}\mu_2-2r^2_{11}\mu_1}{r^2_{21}\Sigma_{11}-2r_{11}r_{21}\Sigma_{12}+\Sigma_{22}r^2_{11}}.
$$

\noindent
\beqlb\label{4-5}
 x^*=\left\{
\begin{array}{ll}
\displaystyle \frac{2r_{11}r_{22}\mu_2-2r^{2}_{22}\mu_1}{\Sigma_{22}r^{2}_{12}-2r_{11}r_{22}\Sigma_{12}+r^2_{22}\Sigma_{11}},
&\;\textrm{if}\; \gamma_2(x_2,Y_0(x))> 0; \\
\infty,&\;\textrm{otherwise},
\end{array}
\right.
 \eeqlb
 where \beqlb\label{4-60}
 Y_0(x)=\frac{-(\mu_2+\Sigma_{22}x)-\sqrt{D_1(x)}}{\Sigma_{22}}.
 \eeqlb
\end{remark}

\subsection{Tail Asymptotics for Marginal Distributions}
In this subsection, we study the tail behaviour of the marginal stationary survival distributions $\P\{Z_i\geq x\}$, $i=1,2$. Here, we only provide detailed analysis for the case that $\P(Z_1\geq x)$. The other case $\P(Z_2\geq x)$ can be similarly studied.

Since $L_1$ only increases at times $t$ for which $\tilde{Z}_1(t)=0$,  from  Corollary \ref{cor-2}, we get that
\beqlb\label{4-48}
    \Phi_1(x) = \Phi_1(0,0)=\frac{\mu_2(r_{12}+\mu_1u_2)-\mu_1(\mu_2u_2+r_{22})}{(r_{11}+\mu_1u_1)(\mu_2u_2+r_{22})-(r_{21}+\mu_2u_1)(u_2\mu_1+r_{12})}\;\text{for any}\;x\in\R.
\eeqlb
Letting $\theta_2=0$ and $\theta_1=x$ in \eqref{2-5} yields
\beqlb\label{4-1}
\Phi(x,0)&&=-\frac{(x r_{11}-u_1\gamma_3(x))\Phi_1(x)+(x r_{12}-u_2\gamma_3(x))\Phi_2(x)}{\gamma_3(x)}\nonumber
\\&&=-\frac{(x r_{11}-u_1\gamma_3(x))\Phi_1(0)+(x r_{12}-u_2\gamma_3(x))\Phi_2(x)}{\gamma_3(x)},
\eeqlb
where \beqlb\label{4-18}
\gamma_3(x)=\mu_1x+\frac{1}{2}\Sigma_{11}x^2=x(\mu_1+\frac{1}{2}\Sigma_{11}x).
\eeqlb
It is obvious that the only non-zero solution of $\gamma_3(x)=0$ is
\beqlb\label{4-2}
x_{\gamma_3}=-\frac{2\mu_1}{\Sigma_{11}}.
\eeqlb
Since $\Sigma_{11}>0$,  the zero $\gamma_3(x)$ is not a pole of $\Phi(x,0)$ if $\mu_1\geq 0$.  In  this case, $\P\{Z_1\geq x\}$  has the same tail asymptotics as $V_2(\cdot)$, referring to the previous section for the analysis based on which of the three candidates: $x^*$, $\tilde{x}$ and $x_2$, would be the dominant singularity. The only difference is the expression for the coefficient, which can be easily determined from Lemma~\ref{thm1} and equation \eqref{4-1}.

Below, we assume that $\mu_1<0$.  From \eqref{4-1}, we know that the asymptotic behaviour of $\Phi(x,0)$ depends on which of the four candidates: $x^*$, $\tilde{x}$, $x_2$ and $x_{\gamma_3}$ would be the dominant singularity. For convenience,  let $z=\min\{x^*,\tilde{x}\}$ and $\alpha_1=\min\{z,x_2,x_{\gamma_3}\}$.

\begin{thm}\label{4-thm2}
Assume that $\mu_1<0$. Then, for the marginal survival distribution $\P\{Z_1\geq z\}$, there are four types of exact tail asymptotics:
\begin{itemize}
\item [(1)] If $\min\{z,x_{\gamma_3}\}< x_2$, $x_{\gamma_3}\neq z$ and $x^*\neq \tilde{x}$; or $\min\{z,x_{\gamma_3}\}< x_2$ and $x_{\gamma_3}= z$ with $Y_0(x^*)=0$\red{;} or $z=x_{\gamma_3}=x_2$; or $z>x_{\gamma_3}=x_2$; or $x_{\gamma_3}>z=x_2$ with $\tilde{x}= x^*$, then,  $\P\big\{Z_1\geq x\big\}$ has an  exponential decay, that is
    \beqnn
    \P\big\{Z_1\geq x\big\}\sim K \exp\{-\alpha_1 x\};
    \eeqnn
\item[(2)] If $\min\{z,x_{\gamma_3}\}< x_2$, $x_{\gamma_3}\neq z$  and  $\tilde{x}=x^*$; or $\min\{z,x_{\gamma_3}\}< x_2$ and $x_{\gamma_3}= z$ with $Y_1(x^*)=0$\red{,} then  $\P\big\{Z_1\geq z\big\}$ has an  exponential decay multiplied by a factor of $x$, that is
    \beqnn
    \P\big\{Z_1\geq x\big\}\sim K x\exp\{-\alpha_1 x\};
    \eeqnn
\item[(3)] If $x_{\gamma_3}>z=x_2$ with $\tilde{x}\neq x^*$, then  $\P\big\{Z_1\geq x\big\}$ has an  exponential decay multiplied by a factor of $x^{-\frac{1}{2}}$, that is
    \beqnn
    \P\big\{Z_1\geq x\big\}\sim K x^{-\frac{1}{2}}\exp\{-\alpha_1 x\};
    \eeqnn
\item[(4)] If $x_2<\min\{z,x_{\gamma_3}\}$, then  $\P\big\{Z_1\geq z\big\}$ has an  exponential decay multiplied by a factor of $x^{-\frac{3}{2}}$, that is
    \beqnn
    \P\big\{Z_1\geq x\big\}\sim K x^{-\frac{3}{2}}\exp\{-\alpha_1 x\}.
    \eeqnn
\end{itemize}

\end{thm}

A detailed proof to Theorem 4.3 is presented in the appendix.

\begin{remark}
The constant $K$ in Theorem~\ref{4-thm2} can be calculated based on Lemma~\ref{thm1} and \eqref{4-1}.
\end{remark}

\subsection{Tail Asymptotics for Joint Distribution along a Direction}
In this subsection, we aim to obtain tail asymptotic properties of the joint stationary distribution $\pi$ along a direction via the kernel method and Theorem \ref{TLT}.  Here, we should point out that the main reason to restrict ourselves to along a direction is because that Theorem~\ref{TLT} only applies to univariate functions. A different method will be used to study tail asymptotics in the joint distribution in Section~\ref{sec:5}. More specifically, we consider, in this section,  the tail behaviour for the distribution of the random variable $<\bar{u},Z>$ along any direction $\bar{u}=(\bar{u}_1,\bar{u}_2)'\in\R_+^2$.

   For the random variable $<\bar{u}, Z>$, we denote the MGF $\Phi_{\bar{u}}(x)$ by
\beqnn
    \Phi_{\bar{u}}(\lambda)=\E_\pi\big[\exp\{\lambda<\bar{u},Z>\}\big].
\eeqnn
It is obvious that
\beqnn
\Phi_{\bar{u}}(\lambda)=\Phi(\lambda \bar{u}):=\Phi(\lambda \bar{u}_1,\lambda \bar{u}_2).
\eeqnn
It follows from the BAR \eqref{2-5} that
\beqlb\label{4-56}
\Phi_{\bar{u}}(\lambda)=-\frac{\Phi_1(\lambda \bar{u}_2)\big(\gamma_1(\lambda \bar{u})-u_1\Psi_X(\lambda \bar{u})\big)+\big(\gamma_2(\lambda \bar{u})-u_2\Psi_X(\lambda \bar{u})\big)\Phi_2(\lambda \bar{u}_1)}{\Psi_X(\lambda \bar{u})}.
\eeqlb
From \eqref{4-56}, we know that asymptotic behaviour of $\Phi_{\bar{u}}(\lambda)$ depends on that of $\Phi_i(\lambda \bar{u}), i=1,2,$ and $\Psi_X(\lambda \bar{u})$.
For convenience, let
\beqnn
\gamma(\lambda)=\Psi_X (\lambda \bar{u})=\lambda<\bar{u},\mu>+\frac{\lambda^2}{2}<\bar{u}, \Sigma \bar{u}>.
\eeqnn
Therefore, the only non-zero solution of $\gamma(\lambda)=0$ is
\beqlb\label{5-57}
x_{\gamma}=-\frac{2<\bar{u},\mu>}{<\bar{u},\Sigma \bar{u}>}.
\eeqlb

If $x_{\gamma}<0$, then the zero $\gamma(\lambda)$ is not a pole of $\Phi(x,0)$.  In  this case, the tail behaviour of the distribution for $<\bar{u}, Z>$  is completely determined by that for $V_i$, which in turn can be obtained directly from Theorem~\ref{thm-a2}. The only difference is the expression for the coefficient, which can be easily obtained from Theorem~\ref{thm-a2} and equation \eqref{4-56}.

Therefore, in the rest of this subsection, we assume that $x_{\gamma}>0$. Let $\beta_i=C_p(\Phi_i(\lambda \bar{u}))$, $i=1,2$.  Then, $\beta_1=\frac{\alpha_{1}}{\bar{u}_2}$ and $\beta_2=\frac{\alpha_{2}}{\bar{u}_1}$.  Below,
without loss of generality, we assume that $\beta_1>\beta_2$. The case of $\beta_1 \leq \beta_2$ can be similarly discussed.  For convenience, let $z_0=\min\{\frac{x^*}{\bar{u}_1},\frac{\tilde{x}}{\bar{u}_1}\}$ and $\beta=\min\{\beta_2,x_{\gamma}\}$.

\begin{thm}\label{thm-2}
Assume that $x_{\gamma}>0$ and $\beta_1>\beta_2$. Then, for the random variable $<\bar{u}, Z>$, we have the following asymptotic behaviour:
\begin{itemize}
\item[(1)] If $\min\{z_0,x_{\gamma}\}< \frac{x_2}{\bar{u}_1}$ with $\tilde{x}=x^*$ and $z_0<x_{\gamma}$; or  $\min\{z_0,x_{\gamma}\}< \frac{x_2}{\bar{u}_1}$ with $z_0=x_{\gamma}$  and  $Y_1(x^*)=x_\gamma \bar{u}_2$, then  $\P\big\{<\bar{u},Z>\geq \lambda\big\}$ has an  exponential decay multiplied by a factor of $\lambda$, that is
    \beqnn
    \P\big\{<\bar{u},Z>\geq \lambda\big\}\sim K \lambda\exp\{-\beta\lambda\}\red{.}
    \eeqnn
\item[(2)] If  $\min\{z_0,x_{\gamma}\}< \frac{x_2}{\bar{u}_1}$ with $\tilde{x}=x^*$ and $z_0> x_{\gamma}$; or $\min\{z_0,x_{\gamma}\}< \frac{x_2}{\bar{u}_1}$ with $\tilde{x}\neq x^*$; or  $\min\{z_0,x_{\gamma}\}< \frac{x_2}{\bar{u}_1}$ with $x_{\gamma}= z_0$ and  $Y_0(x^*)=x_\gamma \bar{u}_2$; or  $\min\{z_0,x_{\gamma}\}=\frac{x_2}{\bar{u}_1}$ with $z_0>x_{\gamma}=\frac{x_2}{\bar{u}_1}$, then
    \beqnn
 \P\{<\bar{u},Z>>\lambda\}\sim K  \exp\{-\beta\lambda\}.
\eeqnn
\item[(3)] If $\min\{z_0,x_{\gamma}\}=\frac{x_2}{\bar{u}_1}$ with $x_{\gamma}>z_0=\frac{x_2}{\bar{u}_1}$ and $\tilde{x}\neq x^*$, then
\beqnn
 \P\{<\bar{u},Z>>\lambda\}\sim K\lambda^{-\frac{1}{2}} \exp\{-\beta\lambda\}.
\eeqnn
\item[(4)] If $\frac{x_2}{\bar{u}_1}<\min\{z_0,x_{\gamma}\}$, then
\beqnn
\P\{<\bar{u},Z>>\lambda\}\sim K \lambda^{-\frac{3}{2}}\exp\{-\beta\lambda\}.
\eeqnn
\end{itemize}
\end{thm}

A detailed proof to Theorem 4.3 is presented in the appendix.

\begin{remark} In Theorem~\ref{thm1-2}, all results are obtained under the assumption $\beta_1>\beta_2$. Results when $\beta_1 \leq \beta_2$ can be similarly obtained.  Here, we provide a brief discussion.
\item[(1)]  {\bf Case  $\beta_1<\beta_2$}: By symmetry, this is the case of interchanging $\beta_1$ and $\beta_2$ in Theorem~\ref{thm1-2}.
\item[(2)] { \bf Case $\beta_1=\beta_2$}: From the proof to Theorem~\ref{thm-2}, we obtain that
\beqlb\label{4-16}
\lim_{\lambda\to\beta_{i+1}}(\beta_{i+1}-\lambda)^{\bar{\mu}_{i+1}}\Phi_{i+1}(\lambda \bar{u}_{2-i})=B(\beta_{i+1},\bar{\mu}_{2-i}),\;i=0,1.
\eeqlb
For convenience, let
\beqnn
\bar{\mu}=\max\{\bar{\mu}_1,\bar{\mu}_2\}.
\eeqnn
From \eqref{4-56}, we get that the tail behaviour of $\Phi_{\bar{u}}(\lambda)$ is determined by $\Phi_i(\lambda)$ which corresponds to $\bar{\mu}$.
From \eqref{4-56} and \eqref{4-16}, we get that
\beqlb\label{4-19}
\lim_{\lambda\to \beta } (\beta-\lambda)^{\bar{\mu}}\Phi_{\bar{u}}(\lambda)=B(\beta,\bar{\mu}).
\eeqlb
Based on \eqref{4-19} and the results for the cases $\beta_1\neq \beta_2$, we can get the results directly.
\end{remark}

\section{Extreme Value Distribution of Joint Distribution} \label{sec:5}

In the previous section, we obtained exact tail asymptotics for the boundary stationary distributions, the marginal distributions and the joint stationary distribution along a direction $\bar{u}\in\R_+^2$ via the kernel method and the Tauberian-like Theorem~\ref{TLT}, which applies only to univariate distributions. Therefore, it cannot be used for tail behaviour of the joint stationary distribution $\pi(x):=F(x)$, which is the focus of this section.

It is well known that the multivariate Gaussian vector with the correlation coefficient being less than $1$ is asymptotic independent. On the other hand, we note that in the interior of the first quadrant, the sticky Brownian motion $Z$ behaves like the Brownian motion $X=(X_1,X_2)'$. Hence, it is expected that $Z$ is also asymptotically independence.  In this section, we prove this fact and study the extreme value of $F(\cdot)$.  In the rest of this paper, we assume that the correlation coefficient $\rho_{X_1X_2}<1$

To achieve our goal, we first note that in the previous section, we obtained the tail equivalence of the marginal distributions. These results provide us with much information for studying the tail dependence of the stationary distribution. Tail dependence describes the amount of dependence in the upper tail or lower tail of a multivariate distributions. Once we clarify their dependence, we can study the bivariate extreme value distribution of the stationary distribution. The extreme value distribution is very useful since from a sample of vectors of maximum, one can make inferences about the upper tail of the stationary distribution using multivariate extreme value theory and copula.

Before we state our main result of this section, we first introduce the domain of attraction of some extreme value distribution function $G(\cdot)$ and asymptotic independence.

\begin{definition}[Domain of Attraction] Assume that $\big\{X_n=(X_n^{(1)},\ldots,X_n^{(d)})'\big\}$ are i.i.d. multivariate random vectors with common distribution $\tilde{F}(\cdot)$ and the marginal distributions $\tilde{F}_i(\cdot)$, $i=1,\ldots,d$. If there exist normalizing constants $a_n^{(i)}>0$ and $b_n^{(i)}\in\R$, $1\leq i\leq d$, $n\geq 1$ such that as $n\to\infty$
\beqnn
\P\Big\{\frac{M_n^{(i)}-b_n^{(i)}}{a_n^{(i)}}\leq x^{(i)},1\leq i\leq d\Big\}&&=\tilde{F}^n\Big(a_n^{(1)}x^{(1)}+b_n^{(1)},\ldots,a_n^{(d)}x^{(d)}+b_n^{(d)} \Big)
\\&&\to G(x^{(1)},\ldots,x^{(d)}),
\eeqnn
where $M_n^{(i)}=\bigvee_{k=1}^{n}X_k^{(i)}$ is the componentwise maxima, then we call the distribution function $G(\cdot)$ a multivariate extreme value distribution function, and
$F$ is in the domain of attraction of $G(\cdot)$. We denote this by $\tilde{F}\in D(G)$.
\end{definition}

\begin{definition}[Asymptotic Independence]  Assume that the extreme value distribution function $G(\cdot)$ has the marginal distributions $G_i(\cdot)$, $i=1,\ldots,d$. If
\beqnn
\tilde{F}^n\Big(a_n^{(1)}x^{(1)}+b_n^{(1)},\ldots,a_n^{(d)}x^{(d)}+b_n^{(d)} \Big)
\to G(x^{(1)},\ldots,x^{(d)}) =\prod_{i=1}^{d}G_i(x^{(i)}),
\eeqnn
then we say that $\tilde{F}(\cdot)$ is asymptotic independent.
\end{definition}

Many efforts have been made to estimate tail probabilities based on the multivariate extreme value distribution. For more information, readers may refer to Foug\'eres~\cite{F2004}, Ledford and  Tawn~\cite{LT1997}, and the references therein.  In this section, we apply the extreme value to study exact tail asymptotics of the joint stationary distribution $F(\cdot)$ of $Z$. As usual, in L\'evy-driven queueing networks, below we assume that the reflection matrix $R=I-P^T$, where $P$ is a substochastic matrix with its spectral radius strictly less than one. Now, we state the main result of this section.

\begin{thm}\label{thm2} For the bivariate sticky Brownian motion $Z=(Z_1,Z_2)'$, we have
\beqlb\label{5-36}
\P\big\{Z_1\geq x, Z_2\geq y\big\}/\Big(K x^{\bar{\mu}_1}y^{\bar{\mu}_2}\exp\big\{-(\alpha_1x+\alpha_2y)\big\}\Big)\to 1,\;\text{as}\;(x,y)\to(\infty,\infty),
\eeqlb
where $\alpha_i$ is the decay parameter associated with $Z_i$ and $\bar{\mu}_i\in\{0,1,-\frac{1}{2},-\frac{3}{2}\}$ is the exponent corresponding to $\alpha_i$ given in Theorem~\ref{4-thm2}.
\end{thm}

In order to prove Theorem~\ref{thm2}, we first consider the extreme copula of the joint stationary distribution $F$ for the sticky Brownian motion $Z$. For this purpose, we first study the bivariate extreme value distribution function of $F(\cdot)$. In fact, we have the following theorem.
\begin{thm}\label{4-thm} For the sticky Brownian motion $Z=(Z_1,Z_2)'$ with the stationary distribution function $F$, we have
\beqnn
F^n(a_n(\bar{\mu}_i,\alpha_i) x^{(i)}+b_n(\bar{\mu}_i,\alpha_i),i=1,2 )\to G_1(x^{(1)})G_1(x^{(2)}),\;\text{as}\; n\to\infty,
\eeqnn
where $a_n(\bar{\mu}_i,\alpha_i)$ and $b_n(\bar{\mu}_i,\alpha_i)$ are given below by \eqref{4-34} and \eqref{4-47}, respectively, and $G_1(x)$ is given by \eqref{5-a1}.
\end{thm}

\begin{remark}From Theorem \ref{4-thm}, we can read that

\noindent{\bf (1)}  $F(\cdot)\in D(G)$; and

\noindent{\bf (2)} $F$ is asymptotic independent.
\end{remark}

We first prove Theorem~\ref{4-thm}, which requires the following technical lemmas.
\begin{lem}\label{4-lem-1}
For any univariate continuous   distribution function $\tilde{F}(\cdot)$,  if $1-\tilde{F}(x)\sim x^{-\bar{\mu}}\exp\{-\alpha x\}$ with $\alpha\in\R_+$ and $\bar{\mu}\in\R$, as $x\to\infty$, then we have
\beqlb\label{5-19}
\tilde{F}'(x)\sim \alpha x^{-\bar{\mu}} \exp\{-\alpha x\},\;\text{as}\; x\to\infty,
\eeqlb
and
\beqlb\label{5-20}
\tilde{F}''(x)\sim -\alpha^2 x^{-\bar{\mu}} \exp\{-\alpha x\},\;\text{as}\; x\to\infty.
\eeqlb
\end{lem}

Since both $\tilde{F}(x)$ and $1-x^{-\bar{\mu}}\exp\{-\alpha x\}$ are differentiable, we can apply L'Hospital's Rule to directly have the result.  Here we skip the proof of this lemma.

In order to study bivariate extreme value distribution of $F(x,y)$, we first need to consider the extreme value distribution functions of the marginal functions $F_i(x)$, $i=1,2$. In fact, we have the following result.
\begin{lem}\label{4-lem-3} For the stationary marginal distribution $F_i(x)=\P\{Z_i\leq x\}$, $i=1,2$, of the $i$th component of the stationary vector of the sticky Brownian motion $Z(t)$, we have
\beqlb\label{4-32}
F_i(x)\in D(G_1),\;i=1,2,\;
\eeqlb
i.e., there exit constants $a_n(\alpha_i,\bar{\mu}_i)$ and $b_n(\alpha_i,\bar{\mu}_i)$, which are functions of $\alpha_i$ and $\bar{\mu}_i$, such that as $n \to \infty$,
\beqnn
F^n_i(a_n(\alpha_i,\bar{\mu}_i)x+b_n(\alpha_i,\bar{\mu}_i))\to G_1(x)
\eeqnn
where \beqlb\label{5-a1} G_1(x)=
\exp\{-e^{-x}\}.
\eeqlb
Moreover, the normalizing constants $a_n(\alpha_i,\bar{\mu}_i)$ and $b_n(\alpha_i,\bar{\mu}_i)$ are given below by \eqref{4-34} and \eqref{4-47}, respectively.
\end{lem}

\proof
We only prove the case $i=1$. The other case can be considered in the same fashion.

It follows from Theorem \ref{4-thm2} that we have
\beqlb\label{5-2}
1-F_1(x)\sim K x^{-\bar{\mu}_1}\exp\{-\alpha_1 x\},\text{as}\;x\to\infty,
\eeqlb
where $\alpha_1\in\{x^*,\tilde{x},x_2\}$ and $\bar{\mu}_1\in\{0,\frac{1}{2},\frac{3}{2},-1\}$.
It follows from the asymptotic equivalence \eqref{5-2} and Lemma~\ref{4-lem-1}  that
\beqlb\label{4-40}
\lim_{x\to\infty}\frac{F''_1(x)\big(1-F_1(x)\big)}{\Big(F'_1(x)\Big)^2}=-1.
\eeqlb
Then, it follows from Proposition~1.1 in Resinck~\cite[pp. 40]{R1987} that $F_1\in D(G_1)$, with $G_1(x)$ given by (\ref{5-a1}).

In the following, we identify suitable normalizing constants $a_n(\bar{\mu}_1,\alpha_1)$  and $b_n(\bar{\mu}_1,\alpha_1)$.  Since we do not know the explicit expression of $F_1(\cdot)$, we apply the tail equivalence to reach our goal.

First, since
\beqlb\label{4-33}
\lim_{x\to\infty}\frac{1-F_1(x)}{F_1'(x)}=\frac{1}{\alpha_1},
\eeqlb
according to Proposition~1.1 in \cite[pp. 40]{R1987}, we can choose
\beqlb\label{4-34}
a_n(\bar{\mu}_1,\alpha_2)= \frac{1}{\alpha_1}.
\eeqlb
Next, we find a suitable $b_n(\bar{\mu}_1,\alpha_1)$. Due to Proposition~1.1 in \cite[pp. 40]{R1987}, one choice of $b_n$ is
\beqlb\label{4-36}
1-F_1(b_n)=\frac{1}{n},
\eeqlb
based on which and the tail asymptotic equivalence, we can choose $b_n$ such that
\beqlb\label{4-41}
  -K\big(b_n\big)^{-\bar{\mu}_1}\exp\{-\alpha_1 b_n\}=\frac{1}{n},\;\textrm{for some constant}\; K,
\eeqlb
i.e.,
\beqlb\label{4-23}
\alpha_1 b_n+\bar{\mu}_1 \log(b_n)+\log(K)=\log (n).
\eeqlb
To identify a solution of $b_n$ to (\ref{4-23}), without loss of generality, we assume that $\bar{\mu}_1\neq 0$ below. Then, we have
 \beqnn
 \log \big(n\big)/b_n\to \alpha_1, \textrm{as}\;b_n\to\infty,
 \eeqnn
 that is
 \beqlb\label{4-42}
 \alpha_1 b_n=\log (n)+r_n,
 \eeqlb
 where $r_n=o\Big(\log(n)\Big)$.

 Combing \eqref{4-23} and \eqref{4-42}, we have
 \beqlb\label{4-43}
 \log(n)+r_n+\bar{\mu}_1 \log\Big(\frac{1}{\alpha_1}\log n+\frac{1}{\alpha_1}r_n\Big)+\log\big(K\big)=\log n.
 \eeqlb
 By some calculations, \eqref{4-43} can be rewritten as
 \beqlb\label{4-44}
 r_n=-\bar{\mu}_1 \log\bigg(\frac{1}{\alpha_1}\Big(1+\frac{r_n}{\log n}\Big)\bigg)-\bar{\mu}_1\log\log\big(n\big)-\log K.
 \eeqlb
 Due to \eqref{4-44}, we have
 \beqlb\label{4-45}
 r_n=-\bar{\mu}_1\log\log\big(n\big)+o(1)-\log K.
 \eeqlb
 Combing \eqref{4-42} and \eqref{4-45}, we have
 \beqnn\label{4-46}
\frac{\alpha_1 b_n+\bar{\mu}_1\log\log\big(n\big)- \log\big(n\big)-\log K}{a_n}=\frac{o(1)}{a_n}\to 0.
 \eeqnn
 Hence, we can choose
 \beqlb\label{4-47}
 b_n=\frac{1}{\alpha_1}\Big(\log\big(n\big)-\bar{\mu}_1\log\log\big(n\big)-\log (K)\Big).
 \eeqlb
Finally, it follows from Proposition~1.19  in Resinck \cite[pp. 72]{R1987} and the convergence to types theorem (see Resinck \cite[Propositions 0.2 and 0.3 ]{R1987})  that we can set $a_n(\alpha_1,\bar{\mu}_1)=a_n$ and $b_n(\alpha_1,\bar{\mu}_1)=b_n$.
\qed

The final piece that we need before we can prove Theorem~\ref{4-thm} is a modified version of Proposition~5.27 in Rensick~\cite[pp. 296]{R1987}.
\begin{lem}\label{4-lem-2}
Suppose that $\big\{X_n=(X_n^{(1)},X_n^{(2)})',\;n\in\mathbb{N}\big\}$ are i.i.d. random vectors in $\R^2$ with the common joint continuous distribution $\tilde{F}(\cdot)$, and the marginal distributions $\tilde{F}_i(\cdot)$, $i=1,2$. Moreover, we assume that $\tilde{F}_i(\cdot)$, $i=1,2$ are both in the domain of attraction of some univariate extreme value distribution $\hat{G}_1(\cdot)$, i.e., there exist constants $a_n^{(i)}$ and $b_n^{(i)}$ such that
\beqnn
    \tilde{F}_i\Big(a_n^{(i)}x+b_n^{(i)}\Big)\to \hat{G}_1(x).
\eeqnn
Then, the following are equivalent:
\begin{itemize}
\item[(1)] $\tilde{F}$ is in the domain of attraction of a product measure, that is,
\beqnn
\tilde{F}^n\Big(a_n^{(1)}x^{(1)}+b_n^{(1)},a_n^{(2)}x^{(2)}+b_n^{(2)}\Big)\to \hat{G}_1\big(x^{(1)}\big)\hat{G}_1\big(x^{(2)}\big);
\eeqnn
\item[(2)]
\beqnn
\P\Big(\bigvee_{l=1}^n X_l^{(1)}\leq a_n^{(1)}x^{(1)}+b_n^{(1)}, \bigvee_{l=1}^n X_l^{(2)}\leq a_n^{(2)}x^{(2)}+b_n^{(2)}\Big)\to \hat{G}_1\big(x^{(1)}\big)\hat{G}_1\big(x^{(2)}\big);
\eeqnn
\item[(3)] For large enough $n\in\mathbb{N}$ such that $a_n^{(i)}x^{(i)}+b_n^{(i)}>0$ with $\hat{G}_1(x^{(i)})>0$, $i=1,2$,
\beqnn
\lim_{n\to\infty}n\P\Big(X_{i}^{(1)}>a_n^{(i)}x^{(i)}+b_n^{(i)},i=1,2\Big)=0;
\eeqnn
\item[(4)] With $\lim_{x\to\infty}\tilde{F}_i(x)=1$,
\beqnn
\lim_{t\to\infty}\P\Big(X^{(1)}>t, X^{(2)}>t\Big)/\big(1-\tilde{F}_i(t)\big)\to 0.
\eeqnn
\end{itemize}
\end{lem}

By a slight modification of the proof to Proposition~5.27 in Rensick~\cite[pp. 296]{R1987}, we can prove the above lemma, details of which is omitted here.

\medskip
Now, we are ready to prove Theorem~\ref{4-thm}.
\medskip

\noindent \underline{\proof of Theorem \ref{4-thm}:}
In the following proof, the reflection matrix $R$, the regulator process $L$ and the 2-dimensional Brownian motion $X=(X_1,X_2)'$ are the components in the definition of the SRBM given in (\ref{def1}), the time change process $T$ is defined through (\ref{2-a3}),  and the sticky Brownian motion $Z$ is defined in (\ref{1-4}). Without loss of generality, we assume that $Z(0)=0$. We mainly use the lemma \ref{4-lem-2} to prove this theorem.
Let
\beqnn
\hat{L}(t)=-[R^{-1}X(t)\wedge R^{-1}\mu t].
\eeqnn
Then it follows from Konstantopoulos, Last and Lin \cite{KLL2004} that  for any $\tilde{z}=(\tilde{z}_1,\tilde{z}_2)'\in\R^2_+$
\beqlb
\P\{Z(t)\geq \tilde{z}\}\leq \P\{\hat{Z}(t)\geq \tilde{z}\}
\eeqlb
where $\hat{Z}(t)=\bar{Z}(T(t))$ with
\beqnn
\bar{Z}(t)=X(t)+R\hat{L}(t).
\eeqnn
It follows from \cite[Lemma 2.7]{K2011} that
\beqnn
T^*=\sup_{\omega}\{T(1,\omega)\}\leq 1,\;\text{a.s}
\eeqnn
By the first change of variable formula (see for example \cite[Proposition 10.21]{J1979}, we have that
\beqnn
\hat{Z}(1)=\int_0^{1} d \hat{Z}(s)=\int_0^{T(1)}d\bar{Z}(s)<\int_{0}^{T^*}d\bar{Z}(s)=\bar{Z}(T^*)\;\textrm{a.s.,}
\eeqnn
where the operations are performed component-wise. Hence for any $\tilde{z}=(\tilde{z}_1,\tilde{z}_2)'\in\R^2_+$
\beqlb\label{5-47}
\P(Z(1)\geq \tilde{z}\}\leq \P\{\bar{Z}(T^*)\geq \tilde{z}\}.
\eeqlb
For convenience, let
\beqnn
\bar{F}(\tilde{z})=\P\big\{Z_1\geq \tilde{z}_1, Z_2\geq \tilde{z}_2\}.
\eeqnn
We also note that
\beqlb\label{5-32}
\bar{F}(\tilde{z})=\lim_{t\to\infty}\P\{Z(t)>\tilde{z}\}=\inf_{t\to\infty}\P\{Z(t)>\tilde{z}\}\leq\P\{Z(1)>\tilde{z}\}.
\eeqlb
From \eqref{5-47} and \eqref{5-32}, we get that
\beqlb\label{5-4}
\bar{F}(\tilde{z})\leq  \P\{\bar{Z}(T^*)\geq \tilde{z}\}
\eeqlb

 Hence, for any $\tilde{z}=(\tilde{z}_1,\tilde{z}_2)\in\R^2_+$ and $t\in\R_+$,
\beqlb\label{4-38}
\P\{Z(T^*)\geq \tilde{z}\}&&\leq \P\{X(T^*)-\mu T^* \geq\tilde{z}\}.
\eeqlb
It is obvious  that $X(T^*)-\mu T^*$ is a Gaussian vector with the correlation coefficient being less than 1.

From \eqref{4-38}, we have that for large enough $z\in\R_+$
\beqlb\label{4-49}
\limsup_{z\to\infty}\frac{\bar{F}(z,z)}{\bar{F}_1(z)}\leq \limsup_{z\to\infty}\frac{ \P\{X(T^*)-\mu T^*\geq (z,z)'\}}{\bar{F}_1(z)}.
\eeqlb

 On the other hand, it is well-known that for any bivariate Gaussian vector with the correlation coefficient being less than one is independent asymptotic.  Hence
\beqlb\label{4-38-a}
&&\limsup_{z\to\infty}\frac{\P\{X_1(T^*)-\mu_1 T^*\geq z, X_2(T^*)-\mu_2 T^*\geq z\}}{\P\{Z_1\geq z\}}\nonumber
\\
&&\hspace{2cm}=\limsup_{z\to\infty}\frac{\P\{X_1(T^*)-\mu_1 T^*\geq z, X_2(T^*)-\mu_2 T^*\geq z\}}{\P\{X_1(T^*)-\mu_1 T^*\geq z\}}\frac{\P\{X_1(T^*)-\mu_1 T^*\geq z\}}{\P\{Z_1\geq z\}}\nonumber
\\&&\hspace{2cm}\leq \limsup_{z\to\infty} \frac{\P\{X_1(T^*)-\mu_1 T^*\geq z, X_2(T^*)-\mu_2 T^*\geq z\}}{\P\{X_1(T^*)-\mu_1 T^*\geq z\}}\nonumber
\\&&\hspace{2cm}\leq \limsup_{z\to\infty} \frac{\P\{X_1(T^*)-\mu_1 T^*\geq z, X_2(T^*)-\mu_2 T^*\geq z\}}{\P\{X_1(T^*)-\mu_1 T^*\geq z\}}=0,
\eeqlb
where the first inequality is obtained by using
\beqnn
\P\{X_1(T^*)-\mu_1 T^*\geq z\}/\P\{Z_1\geq z\}\to 0,\;\text{as}\; z\to\infty.
\eeqnn
From above arguments, we get that
\beqlb\label{4-39}
\lim_{z\to\infty}\frac{\bar{F}(z,z)}{\bar{F}_1(z)}=0.
\eeqlb
The proof to the theorem follows now from  \eqref{4-39} and Lemmas \ref{4-lem-3} and \ref{4-lem-2}.
\qed

Finally, we prove Theorem~\ref{thm2}.
\medskip

\noindent \underline{\proof of Theorem \ref{thm2}:}  To prove this theorem, we first introduce a transformation.
Let $\bar{X}=(\bar{X}_1,\bar{X}_2)'$ be a random vector with joint distribution $\tilde{F}(x,y)$  and marginal distributions $\tilde{F}_i(x)$. Then, we make the following transformation:
\beqlb\label{5-5}
X^*_i=\frac{-1}{\log \big(\tilde{F}_i(\bar{X}_i)\big)},\;\text{for}\;i=1,2.
\eeqlb
By the transformation \eqref{5-5}, we transform each marginal $\bar{X}_i$ of a random vector $\bar{X}$  to a unit Fr\'echet variable $X_i^*$, that is,
\beqnn
\P\{X^*_i<x\}=\exp\{-\frac{1}{x}\}\;\text{for}\;x\in\R_+.
\eeqnn
Hence, for the bivariate extreme value distribution $G(x,y)$
\beqlb\label{5-10}
G^*(x^*,y^*)=G\Big(\frac{-1}{\log \big(G_1(x)\big)},\frac{-1}{\log \big(G_1(y)\big)}\Big),
\eeqlb
where $G^*(\cdot,\cdot)$ is the joint distribution function  with the common marginal Fr\'echnet distribution $\Phi(x)=\exp\{-x^{-1}\}$.
Furthermore,   for the stationary random vector $Z$, define
\beqlb\label{5-12}
Y_i=\frac{1}{1-F_i(Z_i)}.
\eeqlb
Let $F^*(y_1,y_2)$ be the joint distribution function of $Y=(Y_1,Y_2)'$. Then, it follows from Proposition 5.10 in Resnick \cite{R1987} and Theorem~\ref{4-thm} that
\beqlb\label{5-6}
F^*(x,y)\in D\big(G^*(x,y)\big).
\eeqlb
By \eqref{5-6}, we have that for any $Y=(y_1,y_2)'\in\R_+^2$, as $n\to\infty$,
\beqlb\label{5-7}
(F^*(nY))^n\to G^*(Y) .
\eeqlb
It follows from \eqref{5-7} that
\beqnn
F^*(nY)\sim \big(G^*(Y)\big)^{\frac{1}{n}}.
\eeqnn
By a simple monotonicity argument, we can replace $n$ in the above equation by $t$.
Then we have that as $t\to\infty$,
\beqlb\label{5-3}
F^*(tY)\sim \big(G^*(Y)\big)^{\frac{1}{t}}.
\eeqlb
On the other hand, by Lemma \ref{4-lem-3}, for any $y\in\R_+$,
\beqlb\label{5-28}
F^*_i(ty)\sim \big(G^*_1(y)\big)^{\frac{1}{t}},\;\text{for any}\;i=1,2.
\eeqlb
Combing \eqref{5-3} and \eqref{5-28}, we get that as $t\to\infty$
\beqlb\label{5-33}
F^*(tY)\sim F_1^*(ty_1)\cdot F_2^*(ty_2).
\eeqlb
It is obvious that for any $x\in\R_+$
\beqlb\label{5-35}
\bar{F}^*_i(tx):=1-F_i^*(tx)\to 0\;\text{as}\; t\to \infty.
\eeqlb
From \eqref{5-33} and \eqref{5-35}, we have that for any $(x,y)'\in\R_+^2$
\beqlb\label{5-34}
F^*(tx,ty)+\bar{F}^*_1(tx)+\bar{F}^*_2(ty)\sim F_1^*(tx)\cdot F_2^*(ty)+\bar{F}_1^*(tx)+\bar{F}_2^*(ty),
\eeqlb
since
\beqnn
F^*(tx,ty)\to 1,\; F^*_1(tx)\to 1\; \text{and}\; F_2^*(ty) \to 1,\text{as}\;t\to\infty.
\eeqnn
 \eqref{5-34} is equivalent to
\beqlb\label{5-11}
\P\{Y_1\geq tx,Y_2\geq ty\}\sim \bar{F}^*_1(tx)\cdot  \bar{F}^*_2(ty).
\eeqlb
Hence, we have for any $(x,y)'\in\R^2_+$
\beqlb\label{5-44}
\lim_{t\to \infty}\frac{\bar{F}^*(tx,ty)}{\bar{F}^*_1(tx)\cdot \bar{F}^*_2(ty)}=1.
\eeqlb
To prove our theorem, we need to show
\beqlb\label{5-41}
\lim_{(x,y)'\to(\infty,\infty)'}\frac{\bar{F}^*(x,y)}{\bar{F}^*_1(x)\cdot \bar{F}^*_2(y)}=1.
\eeqlb
Note that
\beqlb\label{5-9}
\bar{F}^*(x,y)=\P\big\{\bar{F}_1^*(Y_1)\geq \bar{F}_1^*(x), \bar{F}_2^*(Y_2)\geq \bar{F}_2^*(y) \big\}.
\eeqlb
Hence, to prove \eqref{5-41}, we only need to show
\beqlb\label{5-45}
\lim_{(u,v)'\to(0,0)'}\frac{\hat{C}(u,v)}{uv}=1,
\eeqlb
where $\hat{C}(\cdot,\cdot)$ is the survival copula of $(Y_1,Y_2)'$.

Note that
\beqlb\label{5-8}
\lim_{z\to 0}\frac{\log(1-z)}{z}=1.
\eeqlb
Near the origin $(0,0)'$, the zero sets  of both $\hat{C}(u,v)$ and $uv$  consist of the lines $u=0$ and $v=0$.
For the line $v=0$, choose, say $\vec{z}=(1,1)'$, then $\vec{z}$ is not tangent to the line $v=0$  at the point $(0,0)'$.   Next, we take the direcitonal derivative along the direction $\vec{z}=(1,1)'$. It follows from \eqref{5-44} and \eqref{5-8} that
\beqlb\label{5-42}
\lim_{(u,v)'\to (0,0)'} \frac{D_{\vec{z}}\hat{C}(u,v)}{D_{\vec{z}}(uv)}=1.
\eeqlb
Similar to \eqref{5-42}, for the line $u=0$, along the direction  $\vec{z}=(1,1)'$, we still have
\beqlb\label{5-43}
\lim_{(u,v)'\to (0,0)'} \frac{D_{\vec{z}}\hat{C}(u,v)}{D_{\vec{z}}(uv)}=1.
\eeqlb
From \eqref{5-42}, \eqref{5-43} and  Lawlor \cite{L2012} that
\beqlb\label{5-46}
\lim_{(u,v)'\to (0,0)'} \frac{\hat{C}(u,v)}{uv}=1.
\eeqlb
Finally, it follows from \eqref{5-12} that for any $(x,y)'\in\R_+^2$,
\beqlb\label{5-13}
\P\{Z_1\geq x,Z_2\geq y\}&&=\P\big\{F_1(Z_1)\geq F_1(x),F_2(Z_2)\geq F_2(y)\big\}\nonumber
\\&&=\P\big\{Y_1\geq \frac{1}{1-F_1(x)},Y_2\geq \frac{1}{1-F_2(y)}\big\}\nonumber
\\&&=F^*\Big(\frac{1}{\bar{F}_1(x)},\frac{1}{\bar{F}_2(y)}\Big).
\eeqlb
Combining \eqref{5-41} and \eqref{5-13}, we get
\beqlb\label{5-15}
\P\{Z_1\geq x,Z_2\geq y\}/\Bigg( \bar{F}^*_1\bigg(\frac{1}{\bar{F}_1(x)}\bigg)\cdot  \bar{F}^*_2\bigg(\frac{1}{\bar{F}_2(y)}\bigg)\Bigg)\to 1,\;\text{as}\; (x,y)'\to(\infty,\infty)'.
\eeqlb
 On the other hand, we get that
\beqlb\label{5-39}
\lim_{x\to 0}\frac{1-\exp\{-x\}}{x}=1.
\eeqlb
By \eqref{5-15} and \eqref{5-39}, we get that
\beqlb\label{5-40}
\P\{Z_1\geq x,Z_2\geq y\}/\Big( \bar{F}_1(x)\cdot \bar{F}_2(y)\Big)\to 1,\;\text{as}\; (x,y)'\to (\infty,\infty)'.
\eeqlb
Finally, it follows from Theorem \ref{4-thm2}  and \eqref{5-40} that
\beqlb\label{5-16}
\bar{F}_i(x)\sim K x^{\bar{\mu}_i}\exp\{-\alpha_i x\},\;i=1,2.
\eeqlb
From above arguments, the theorem is proved. \qed

\section{Final Note} \label{sec:6}

In this work, we studied tail properties of stationary distributions for a two-dimensional sticky Brownian motion \eqref{1-4}, which is a time-changed SRBM \and an extension of SRBM.  Tail asymptotics for stationary distributions of SRBM have attracted a lot of interest recently. For example, for tail asymptotic properties in a boundary distribution or for the joint distribution along a direction (path), Dai and Miyazawa~\cite{DM2011,DM2013} used a geometric method, Dai, Dawson and Zhao~\cite{DDZ} extended the kernel method, Franceschi and Kurkova~\cite{FK2016} employed analytic methods, and Franceschi and Raschel~\cite{FR2017} applied the boundary value problem. In this paper, we made the following contributions:
\begin{itemize}
\item[(1)] A comprehensive study on tail asymptotics for the two-dimensional sticky Brownian motion. The results for the boundary stationary distributions obtained in this paper in parallel to those reported in \cite{DDZ}. In addition, we also considered exact tail asymptotics for the marginal distributions and the joint distribution along an arbitrary direction.
\item[(2)] A study of exact tail asymptotics for the joint distribution using a different method, a combination of the extreme theory and copula.
\end{itemize}

We expect that the kernel method and especially the method combining with the extreme value theory and copula can be extended for studying exact tail asymptotic properties for many other stochastic models.

\noindent\\[0.1mm]

\noindent{\bf Acknowledgments:} This research work was supported by the National Natural Science Foundation of China (No.11361007), the Fostering Project of Dominant Discipline and Talent Team of Shandong Province Higher Education Institutions, and the Natural Sciences and Engineering Research Council (NSERC) of Canada. The authors also thank Dr. Haijun Li for pointing us to references on copula.

\begin{appendix}
\section{Appendix: Technical Proofs}
\subsection{Proof of Theorem \ref{4-thm2}}

\noindent \underline{\proof of Theorem \ref{4-thm2}:}
 Under the assumption of $\mu_1<0$, $x_{\gamma_3}$ might be a pole of $\Phi(x,0)$. Hence, we need to consider the relationship between $x^*$, $\tilde{x}$, $x_2$ and $x_{\gamma_3}$, and to consider the following cases:
\medskip

  \noindent{{\bf Case 1}: $\min\{z,x_{\gamma_3}\}< x_2$ and $x_{\gamma_3}\neq z$}.

  \noindent{{\bf Subcase 1-1}: $\tilde{x}=x^*$ and $z<x_{\gamma_3}$}. In this case, we have  $\alpha_{1}=z$. It follows from Lemma \ref{thm1} and  \eqref{4-1} that
   \beqnn
 \lim_{x\to x^*}(x^*-x)^2\Phi(x,0)&&=\lim_{x\to x^*}\frac{(2r_{12}+u_2\Sigma_{11}(x_{\gamma_3}-x))\Phi_2(x)}{\Sigma_{11}(x_{\gamma_3}-x)}(x^*-x)^2
 \\&&=\frac{(2r_{12}+u_2\Sigma_{11}(x_{\gamma_3}-x^*))A_4(x^*)}{\Sigma_{11}(x_{\gamma_3}-x^*)}.
 \eeqnn
We can check that $\Phi(x,0)$ satisfies conditions in Theorem~\ref{TLT}. Hence
  by Theorem~\ref{TLT},
  \beqlb\label{4-53}
  \P\{Z_1\geq x)\sim  K x \exp\{-\alpha_{1}x\}.
  \eeqlb

\noindent{{\bf Subcase 1-2}:  $\tilde{x}=x^*$ and   $z>x_{\gamma_3}$}. In this case, we have $\alpha_{1}=x_{\gamma_3}$, which implies that $\Phi_2(x)$ is analytic at $ x_{\gamma_3}$.
It follows from Lemma \ref{thm1} and  \eqref{4-1} that
   \beqlb\label{4-51}
 \lim_{x\to x_{\gamma_3}}(x_{\gamma_3}-x)\Phi(x,0)&&=\frac{2r_{11}\Phi_1(0)+2r_{12}\Phi_2(x_{\gamma_3})}{\Sigma_{11}}.
 \eeqlb
Therefore, by Theorem~\ref{TLT} we have
\beqnn
\P\{Z_1\geq x)\sim  K \exp\{-\alpha_{1}x\}.
\eeqnn

  \noindent{{\bf Subcase 1-3}:
$\tilde{x}\neq x^*$}. In this case $\alpha_1=\min\{z,x_{\gamma_3}\}$.  Moreover, $\tilde{x}$, $x^*$ and $x_{\gamma_3}$ are all different.
  So, $\alpha_1$ is a single pole of $\Phi(x,0)$.  Hence,
it follows from Lemma~\ref{thm1}, \eqref{4-1} and \eqref{4-51} that
   \beqlb\label{4-52}
 \lim_{x\to \alpha_{1}}(\alpha_1-x)\Phi(x,0)&&=K(\alpha_{1}),
 \eeqlb
where $K(\alpha_1)$is a constant depending on $\alpha_{1}$.  It then follows from Theorem~\ref{TLT} that
  \beqlb\label{4-10}
  \P\{Z_1\geq x)\sim K \exp\{-\alpha_{1}x\}.
  \eeqlb
\medskip

\noindent{{\bf Case 2}: $\min\{z,x_{\gamma_3}\}< x_2$ and $x_{\gamma_3}= z$}.  In this case, $\alpha_1=x_{\gamma_3}$. We show that we must have
   \beqlb\label{4-11}
   x_{\gamma_3}\neq \tilde{x}.
   \eeqlb
Otherwise, if
\beqnn
 \tilde{x}=x_{\gamma_3},
 \eeqnn
$\gamma_3(x)=\Psi_{X}(x,0)$ implies that
    \beqlb\label{4-12}
  Y_0(\tilde{x})=0\;\text{or}\; Y_1(\tilde{x})=0.
    \eeqlb
 On the other hand, it follows from the proof of Lemma \ref{4-lem-5}  that in this case
  $Y_0(\tilde{x})>0$, which contradicts to \eqref{4-12}.  Hence, we must have \eqref{4-11}, which is equivalent to
 \beqlb\label{4-3}
 x^*=x_{\gamma_3}.
 \eeqlb
Once again, since $\gamma_3(x)=\Psi_{X}(x,0)$, \eqref{4-3} implies that
  \beqlb\label{4-12a}
  Y_0(x^*)=0\;\text{or}\; Y_1(x^*)=0.
    \eeqlb

\noindent{\bf Subcase 2-1:} Assume that $Y_0(x^*)=0$. Recall the definition of $x^*$, which satisfies
 \beqlb\label{4-8}
 \gamma_2(x^*,Y_0(x^*))=0.
 \eeqlb
Therefore, we have
 \beqlb\label{4-7}
 \gamma_2(x^*,0)=x^*r_{12}=0
 \eeqlb
since $Y_0(x^*)=0$.
 Since $x^*>0$, we must have $r_{12}=0$. It follows from the above arguments that equation \eqref{4-1} is equivalent to
 \beqlb\label{4-14}
 \Phi(x,0)=-\frac{x r_{11}\Phi_1(0)}{\gamma_3(x)}+u_1\Phi_1(0)+u_2\Phi_2(x).
 \eeqlb

It follows from \eqref{thm1-1}, \eqref{4-18} and \eqref{4-14} that  $\Phi(x,0)$ has a single pole $x^*$. Hence,
\beqnn
\lim_{x\to x_{\gamma_3}}(x_{\gamma_3}-x)\Phi(x,0)=u_2A_1(x^*)+\frac{2r_{11}\Phi_1(0)}{\Sigma_{11}}.
\eeqnn

It is easy to check that $\Phi(x,0)$ satisfies conditions in Theorem \ref{TLT}.
By Theorem \ref{TLT}, we conclude that
\beqlb\label{4-15}
\P\big(Z_1\geq x\big)\sim K \exp\{-\alpha_{1}x\}.
\eeqlb

\noindent{\bf Subcase 2-2:} Assume that $Y_1(x^*)=0$.   In this case
\beqnn
\gamma_2(x^*,0)\neq 0.
\eeqnn
It follows from Lemma \ref{thm1} and \eqref{4-1} that
\beqnn
\lim_{x\to x_{\gamma_3}}(x_{\gamma_3}-x)^2\Phi(x,0)=\lim_{x\to x_{\gamma_3}}(x_{\gamma_3}-x)^2\frac{(2r_{12}+u_2\Sigma_{11}(x_{\gamma_3}-x))\Phi_2(x)}{\Sigma_{11}(x_{\gamma_3}-x)}=\frac{2r_{12}A_1(x^*)}{\Sigma_{11}}.
\eeqnn
By Theorem \ref{TLT}, we conclude that
\beqnn
\P\big(Z_1\geq x\big)\sim K x \exp\{-\alpha_{1}x\}.
\eeqnn
\medskip

\noindent{\bf{Case 3:}} $\min\{z,x_{\gamma_3}\}=x_2$\red{.}  In this case, the following four possibilities exist.

 \noindent{\bf{Subcase 3-1: $z=x_{\gamma_3}=x_2$ }}. In this subcase, we have $\alpha_1=x_2=x_{\gamma_3}$. Similar to Case~2, we only have
 \beqlb\label{4-20}
 x^*=x_{\gamma_3}=x_2.
 \eeqlb
 Moreover,
 \beqnn
 Y_0(x^*)=Y_1(x^*).
 \eeqnn
 Hence, \eqref{4-14} is still valid.  Therefore,
 \beqnn
\lim_{x\to x_{\gamma_3}}(x_{\gamma_3}-x)\Phi(x,0)=\frac{2r_{11}\Phi_1(0)}{\Sigma_{11}}.
\eeqnn
By Theorem~\ref{TLT}, we have
\beqnn
\P\big(Z_1\geq x\big)\sim K \exp\{-\alpha_{1}x\}.
\eeqnn

 \noindent{\bf{Subcase 3-2: $z>x_{\gamma_3}=x_2$ }}.  In this subcase, $\alpha_1=x_{\gamma_3}=x_2$ and $\Phi_1\big(Y_0(x)\big)$ is analytic at $\hat{y}=Y_0(x_2)$. Moreover,
$\gamma_2\big(x_2,Y_0(x_2)\big)\neq 0$.
 From \eqref{4-1} and \eqref{4-18}, we have
 \beqlb\label{4-21}
 \lim_{x\to x_2}(x_2-x)\Phi(x,0)&&=\lim_{x\to x_2}(x_2-x)\frac{(2r_{12}+u_2\Sigma_{11}(x_2-x))\Phi_2(x)}{\Sigma_{11}(x_2-x)}+\frac{2r_{11}\Phi_1(0)}{\Sigma_{11}}\nonumber
 \\&&=\lim_{x\to x_2}\gamma_1\big(x,Y_0(x)\big)\Phi_1\big(Y_0(x)\big)\frac{(2r_{12}+u_2\Sigma_{11}(x_2-x))(x_2-x)}{\Sigma_{11}\gamma_2\big(x,Y_0(x)\big)(x_2-x)}+\frac{2r_{11}\Phi_1(0)}{\Sigma_{11}}\nonumber
 \\&&=\gamma_1\big(x_2,Y_0(x_2)\big)\Phi_1\big(Y_0(x_2)\big)\frac{2r_{12}}{\Sigma_{11}\gamma_2\big(x_2,Y_0(x_2)\big)}+\frac{2r_{11}\Phi_1(0)}{\Sigma_{11}}.
 \eeqlb
By Theorem~\ref{TLT} we have
\beqlb\label{4-22}
 \P\{Z_1\geq x\}\sim K \exp\{-\alpha_{1}x\}.
\eeqlb

 \noindent{\bf {Subcase 3-3: $x_{\gamma_3}>z=x_2$ with $\tilde{x}\neq x^*$}}.  In this subcase, $\alpha_1=x_2$.
 It follows from \eqref{thm1-2} and \eqref{4-1} that
 \beqnn
 \lim_{x\to x_2}\sqrt{x_2-x}\Phi(x,0)&&=\lim_{x\to x_2}\sqrt{x_2-x}\frac{(2r_{12}+u_2\Sigma_{11}(x_{\gamma_3}-x))\Phi_2(x)}{\Sigma_{11}(x_{\gamma_3}-x)}
 \\&&=\frac{(2r_{12}+u_2\Sigma_{11}(x_{\gamma_3}-x_2))A_2(x_2)}{\Sigma_{11}(x_{\gamma_3}-x_2)}.
 \eeqnn
Then by Theorem \ref{TLT}, we have
\beqnn
 \P(Z_1\geq x)\sim K x^{-\frac{1}{2}} \exp\{-\alpha_{1}x\}.
\eeqnn

 \noindent{\bf{Subcase 3-4: $x_{\gamma_3}>z=x_2$ with $\tilde{x}= x^*$}}.  In this subcase, $\alpha_1=x_2=z$.
 From Lemma \ref{thm1} and \eqref{4-1},
  \beqnn
 \lim_{x\to x_2}(x_2-x)\Phi(x,0)&&=\lim_{x\to x_2}\frac{(2r_{12}+u_2\Sigma_{11}(x_{\gamma_3}-x))\Phi_2(x)}{\Sigma_{11}(x_{\gamma_3}-x)}
 \\&&=\frac{(2r_{12}+u_2\Sigma_{11}(x_{\gamma_3}-x_2))A_1(x_2)}{\Sigma_{11}(x_{\gamma_3}-x_2)}.
 \eeqnn
Then by Theorem \ref{TLT}, we have
\beqnn
 \P(Z_1\geq x)\sim K  \exp\{-\alpha_{1}x\}.
\eeqnn
\medskip

 \noindent{\bf{Case 4: $x_2<\min\{z,x_{\gamma_3}\}$ }}.  In this case, $\alpha_{1}=x_2$.
 We first show that in this case
 \beqlb\label{4-26}
 \Phi_2(x_2)<\infty.
 \eeqlb
 In fact, since $x_2<\min\{z,x_{\gamma_3}\}$, we have
 \beqnn
 \gamma_2(x_2, Y_0(x_2))\neq 0,\;\text{and}\; \gamma_1(x_2,Y_0(x_2))\neq 0.
 \eeqnn
Moreover, we get that $\Phi_1(Y_0(x))$ is analytic at $x_2$. From above arguments and Lemma \ref{4-lem-5}, we can get \eqref{4-26}.

On the other hand, it follows from equation \eqref{thm1-3} that
 \beqlb
 \Phi'_2(x)\sim K (x_2-x)^{-\frac{1}{2}},\;\text{as}\;x\to x_2.
 \eeqlb
Hence,
 \beqlb\label{4-a1}
\lim_{x\to x_2} \frac{\int_x^{x_2} \Phi_2'(y)dy}{\int_x^{x_2} K (x_2-u)^{-\frac{1}{2}}du}=1.
 \eeqlb
Therefore,
\beqlb\label{4-54}
\Phi_2(x_2)-\Phi_2(x)\sim K \sqrt{x_2-x}.
\eeqlb

By  \eqref{4-1} and \eqref{4-54}, we get that
\beqnn
\lim_{x\to x_2} \frac{\Phi(x_2,0)-\Phi(x,0)}{K\sqrt{x_2-x}}= 1,
\eeqnn
that is
\beqlb\label{4-27}
\lim_{x\to x_2} K \sqrt{x_2-x} \frac{\Phi(x_2,0)-\Phi(x,0)}{x_2-x}= 1.
\eeqlb
From Dai and Miyazawa~\cite{DM2011}, we get that  $$\frac{\Phi(x_2,0)-\Phi(x,0)}{x_2-x}$$ is the moment generating function of the density function
\beqlb\label{s6-6}
\bar{f}(x)=e^{-x_2x}\int_x^\infty e^{x_2u}f(u)du,
\eeqlb
 where $f(x)$  is the density function of the marginal distribution $\P(Z_1\leq x)$.  Therefore, from Theorem \ref{TLT} and \eqref{4-27}, we have
 \beqlb\label{s6-5}
 \bar{f}(z)\sim K z^{-\frac{1}{2}}e^{-x_2 z}.
 \eeqlb
From \eqref{s6-5} and \eqref{s6-6}, we obtain that
 \beqlb\label{s6-9}
 \int_x^\infty e^{x_2u}f(u)du\sim K x^{-\frac{1}{2}}.
 \eeqlb
 Then, from \eqref{s6-9}, we get
 \beqlb\label{s6-7}
 f(x)\sim K e^{-x_2 x}x^{-\frac{3}{2}}.
 \eeqlb
Hence,
\beqlb
\lim_{x\to\infty}\frac{\int_x^\infty f(u)du}{\int_x^\infty K e^{-x_2 u}u^{-\frac{3}{2}}du}=1,
\eeqlb
and therefore
\beqlb\label{4-84}
 \P(Z_1\geq x)\sim K  x^{-\frac{3}{2}}\exp\{-\alpha_{1}x\}.
\eeqlb
From above arguments, the proof to the theorem is now complete. \qed

\subsection{Proof of Theorem \ref{thm-2}}
\noindent \underline{\proof of Theorem \ref{thm-2}:}
To prove this theorem, we need to consider the relationship between $\beta_2=\min\{\frac{x_2}{\bar{u}_1}, \frac{x^*}{\bar{u}_1}, \frac{\tilde{x}}{\bar{u}_1} \}$ and $x_{\gamma}$. We first point out that, since  $\beta_1>\beta_2$, $\Phi_1(\lambda \bar{u}_2)$ is analytic at $\beta_2$.

 Similar to the previous subsection, we consider the following cases:

 \noindent{{\bf Case 1}: $\min\{z_0,x_{\gamma}\}< \frac{x_2}{\bar{u}_1}$ with $x_{\gamma}\neq z_0$}.   It is obvious that $\beta_2=\{z_0,x_\gamma\}$ in this case.  We consider the following subcases respectively.

  \noindent{{\bf Subcase 1-1} $\tilde{x}=x^*$} and $z_0<x_{\gamma}$. In this case, we have $\beta=z_0$.   It follows from the proof of Lemma \ref{4-lem-5} that
  \beqlb\label{4-61}
  \gamma_2(x^*, Y_0(x^*))=0.
  \eeqlb
There are two possibilities:
 \beqnn
 Y_0(x^*)=\frac{x^* \bar{u}_2}{\bar{u}_1}\;\text{or}\; Y_0(x^*)\neq\frac{x^* \bar{u}_2}{\bar{u}_1}.
 \eeqnn

 We first assume that $Y_0(x^*)=\frac{x^* \bar{u}_2}{\bar{u}_1}$.
  From \eqref{4-61} we have
 \beqlb\label{4-57}
 \gamma_2(\lambda^* \bar{u})=0=\lambda^* (\bar{u}_1r_{12}+\bar{u}_2r_{22}),
 \eeqlb
 where $\lambda^*=\frac{x^*}{\bar{u}_1}$.
 Since $\lambda^*>0$, we must have
 \beqlb\label{4-58}
  (\bar{u}_1r_{12}+\bar{u}_2r_{22})=0.
 \eeqlb
 On the other hand,  since $x^*=\tilde{x}$,  we have that $Y_0(x^*)$ is a pole of $\Phi_1(x)$.  Then, from the assumption, we have
 \beqlb\label{4-29}
\beta_1= Y_0(x^*)/\bar{u}_2>\beta_2=\lambda^*,
 \eeqlb
 which is impossible.  Hence, we can not have $Y_0(x^*)=\frac{x^* \bar{u}_2}{\bar{u}_1}$ and must have
 \beqnn
 Y_0(x^*)\neq\frac{x^* \bar{u}_2}{\bar{u}_1}.
 \eeqnn
 Therefore,
 \beqnn
 \gamma_2(\lambda^* \bar{u})\neq 0.
 \eeqnn
Hence, from Lemma \ref{thm1},  \eqref{4-56} and above arguments, we get that
   \beqnn
 \lim_{\lambda\to \beta}(\beta-\lambda)^2\Phi_{\bar{u}}(\lambda)&&=\lim_{\lambda\to \beta}\frac{\big(\gamma_2(\lambda \bar{u})-u_2\gamma(\lambda)\big)\Phi_2(\lambda \bar{u}_1)}{\gamma(\lambda)}(\beta-\lambda)^2
 \\&&=\frac{\big(\gamma_2(\beta \bar{u})-u_2\gamma(\beta)\big)A_4(x^*)}{\gamma(\beta)}.
 \eeqnn
Moreover, from \eqref{4-56}, it is obvious  that $\Phi_{\bar{u}}(\lambda)$ satisfies  Theorem \ref{TLT}. Hence
  \beqlb\label{5-58}
  \P\{<\bar{u},Z>\geq \lambda)\sim K \lambda \exp\{-\beta \lambda\}.
  \eeqlb

\noindent{{\bf Subcase 1-2} $\tilde{x}=x^*$} and $z_0> x_{\gamma}$.  Here we have $\beta=x_\gamma$.
In such a case, $\Phi_2(\lambda \bar{u}_1)$ is analytic at $ x_{\gamma}$.
It follows from Lemma \ref{thm1} and  \eqref{4-56} that
\beqlb\label{4-59}
 \lim_{\lambda\to \beta}(\beta-\lambda)\Phi_{\bar{u}}(\lambda)&&=2\frac{\gamma_2(\beta \bar{u})\Phi_2(\beta \bar{u}_1)+\gamma_1(\beta \bar{u})\Phi_1(\beta \bar{u}_2)}{ \beta<\bar{u},\Sigma \bar{u}>}.
 \eeqlb
Since $\beta$ is a single pole of $\Phi_{\bar{u}}(\lambda)$, we can easily get that $\Phi_{\bar{u}}(\lambda)$ satisfies the conditions in Theorem \ref{TLT}. Hence
\beqnn
\P\{<\bar{u},Z>\geq \lambda)\sim K \exp\{-\beta \lambda\}.
\eeqnn

\noindent{{\bf Subcase 1-3}: $\tilde{x}\neq x^*$}.  In this case, $\beta=\min\{z_0, x_\gamma\}$. Moreover, from Lemma \ref{thm1} and \eqref{4-56}, we get that $\beta$ is a single pole of $\Phi_{\bar{u}}(\lambda)$. Hence,
\beqlb\label{4-68}
 \lim_{\lambda\to \beta}(\beta-\lambda)\Phi_{\bar{u}}(\lambda)&&=K(\beta),
 \eeqlb
where $K(\beta)$ is a constant depending on $\beta$.
By Theorem \ref{TLT},
  \beqlb\label{4-69}
  \P\{<\bar{u},Z>\geq \lambda)\sim K \exp\{-\beta \lambda\}.
  \eeqlb

\noindent{{\bf Case 2}: $\min\{z_0,x_{\gamma}\}< \frac{x_2}{\bar{u}_1}$ and $x_{\gamma}= z_0$}.  We first point out that in such a case we cannot  have
   \beqlb\label{4-70}
   x_{\gamma}= \tilde{x}/\bar{u}_1=x^*/\bar{u}_1.
   \eeqlb
In fact since $ \gamma(x_\gamma)=\Psi_{X}(x_\gamma \bar{u}_1,x_\gamma \bar{u}_2)=0$, we have that
    \beqlb\label{4-71}
  Y_0(x_\gamma \bar{u}_1)=x_\gamma \bar{u}_2\;\text{or}\; Y_1(x_\gamma \bar{u}_1)=x_\gamma \bar{u}_2.
    \eeqlb
On the other hand, it follows from the proof of Lemma \ref{4-lem-5} that
\beqnn
\gamma_2(\bar{u}_1 x_\gamma, Y_0(\bar{u}_1 x_\gamma))=0,
\eeqnn
that is
\beqlb\label{4-63}
Y_0(x_\lambda \bar{u}_1)=-\frac{x_\lambda \bar{u}_1 r_{12}}{r_{22}}=-\frac{\tilde{x}r_{12}}{r_{22}}.
\eeqlb
From Lemma \ref{4-lem-5} again, in this case, we should have
\beqlb\label{4-64}
Y_0(\tilde{x})>0.
\eeqlb
Combing \eqref{4-63} and \eqref{4-64} gives
\beqlb\label{4-65}
 r_{12}<0.
\eeqlb
since $\bar{u}_1>0$, $r_{22}>0$ and $\tilde{x}>0$.
From Lemma \ref{4-lem7},
\beqlb
X_1(-\frac{x_\lambda \bar{u}_1 r_{12}}{r_{22}})=x_\gamma \bar{u}_1.
\eeqlb
From Lemma \ref{4-lem-5} and \eqref{4-63},
\beqnn
\gamma_1(X_0(-\frac{x_\lambda \bar{u}_1 r_{12}}{r_{22}}), -\frac{x_\lambda \bar{u}_1 r_{12}}{r_{22}})=0,
\eeqnn
that is
\beqnn
X_0(-\frac{x_\lambda \bar{u}_1 r_{12}}{r_{22}})=\frac{x_\lambda \bar{u}_1 r_{12}r_{21}}{r_{11}r_{22}}>0.
\eeqnn
From Remark \ref{cor-1}, we get that
\beqlb
\frac{x_\lambda \bar{u}_1 r_{12}r_{21}}{r_{11}r_{22}}<x_\gamma \bar{u}_1,
\eeqlb
which contradicts to \eqref{4-70}.  Without loss of generality, in this case we assume that
\beqlb\label{4-81}
\frac{x^*}{\bar{u}_1}=x_\gamma.
\eeqlb

In this case, we have $\beta=x_\gamma$. From \eqref{4-71}, there are two subcases in this case.

 \noindent{{\bf Subcase 2-1: $Y_0(x^*)=x_\gamma \bar{u}_2$.}} It follows from Lemma \ref{4-lem-5} that
 \beqlb\label{4-72}
 \gamma_2(x^*,Y_0(x^*))=\gamma_2(x_\gamma \bar{u}_1,x_\gamma \bar{u}_2)=0.
 \eeqlb
 From \eqref{4-72}, we can get \eqref{4-58}.  Moreover, we have
 \beqlb\label{4-28}
 \gamma_1(\beta \bar{u})\neq 0.
 \eeqlb
 If \eqref{4-28} is not true, then we have
 \beqlb\label{4-67}
 \bar{u}_1r_{11}+\bar{u}_2r_{21}=0.
 \eeqlb
 Combing \eqref{4-58} and \eqref{4-67}, we get
 \beqnn
 \bar{u}_2(r_{22}r_{11}-r_{12}r_{21})=0,
 \eeqnn
 which contradicts to Remark \ref{cor-1} and the assumption that $\bar{u}>0$.

  Hence
 it follows from \eqref{thm1-1} and \eqref{4-56} that
 \beqlb\label{4-73}
\lim_{\lambda\to \beta}(\beta-\lambda) \Phi_{\bar{u}}(\lambda)=\frac{2\Phi_1(\beta \bar{u}_2)\gamma_1(\beta \bar{u}_1,\beta \bar{u}_2)}{\beta<\bar{u},\Sigma \bar{u}>}.
 \eeqlb
 From \eqref{4-56} and  \eqref{4-58}, we can easily get that $\Phi_{\bar{u}}(\lambda)$ satisfies Theorem \ref{TLT}. Hence,
 \beqlb\label{4-74}
\P\big(<\bar{u},Z>\geq\lambda\big)\sim K \exp\{-\beta \lambda\}.
\eeqlb

\noindent{{\bf Subcase 2-2: $Y_1(x^*)=x_\gamma \bar{u}_2$.}}  Note that for fixed $x$, $\gamma_2(x,y)$ is increasing in $y$. On the other hand,  for $x^*\in(0,x_2)$, $Y_0(x^*)<Y_1(x^*)$. Therefore, we have
\beqnn
\gamma_2\big(x^*,Y_1(x^*)\big)>\gamma_2(x^*,Y_0(x^*))=0.
 \eeqnn
 From \eqref{4-56},
 \beqlb\label{4-17}
\lim_{\lambda\to \beta}(\beta-\lambda)^2\Phi_{\bar{u}}(\lambda)&&=\lim_{\lambda\to \beta}-(\beta-\lambda)^2\frac{\big(\gamma_{2}(\lambda \bar{u})-u_2\gamma(\lambda)\big)\Phi_2(\lambda \bar{u}_1)}{\gamma(\lambda)}\nonumber
\\&&=\lim_{\lambda\to \beta}2(\beta-\lambda)\frac{\gamma_{2}(\lambda \bar{u})\Phi_2(\lambda \bar{u}_1)}{<\bar{u},\Sigma \bar{u}>}\nonumber
\\&&=\frac{2\gamma_{2}(\beta \bar{u})A_1(\beta \bar{u})}{\beta<\bar{u},\Sigma \bar{u}>}.
 \eeqlb
In this case, $\beta$ is a double pole of $\Phi_{\bar{u}}(\lambda)$. Then, it is easy to verify that $\Phi_{\bar{u}}(\lambda)$ satisfies Theorem \ref{TLT}.
Hence,
\beqlb\label{4-75}
\P\big(<\bar{u},Z>\geq \lambda\big)\sim K \lambda \exp\{-\beta \lambda\}.
\eeqlb

 \noindent{\bf{Case 3: $\min\{z_0,x_{\gamma}\}=\frac{x_2}{\bar{u}_1}$ }}. In this case,  we have $\beta=\frac{x_2}{\bar{u}_1}$ and the following possibilities exist.

 \noindent{\bf{Subcase 3-1: $z_0=x_{\gamma}=\frac{x_2}{\bar{u}_1}$ }}. Similar to Case 2, we cannot  have
 \beqnn
 \frac{x^*}{\bar{u}_1}=x_{\gamma}=\frac{\tilde{x}}{\bar{u}_1}.
 \eeqnn
 Without loss of generality, we assume that
 \beqlb\label{4-76}
 \frac{x^*}{\bar{u}_1}=x_{\gamma}\neq \frac{\tilde{x}}{\bar{u}_1}.
 \eeqlb
 Since $x^*=x_2$,
 \beqlb\label{4-82}
 Y_0(x^*)=Y_1(x^*).
 \eeqlb
On the other hand, we have
\beqlb\label{4-86}
\gamma(x_\gamma)=\Psi(x_\gamma \bar{u}_1,x_\gamma \bar{u}_2)=0.
\eeqlb
Hence, from \eqref{4-82} and \eqref{4-86},
 \beqlb\label{4-83}
 Y_0(x^*)=Y_1(x^*)=x_\gamma \bar{u}_2.
 \eeqlb
 On the other hand, it follows from Lemma \ref{4-lem-5} and above arguments that $Y_0(x^*):=y_0$ is a single pole of $\Phi_1(y)$.  Hence,
 \beqnn
 \beta_1=Y_0(x^*)/\bar{u}_2=x_\gamma=\beta_2,
 \eeqnn
 which contradicts to our assumption that $\beta_1>\beta_2$.  Therefore, we cannot have such a case.

\noindent{\bf{Subcase 3-2 $z_0>x_{\gamma}=\frac{x_2}{\bar{u}_1}$ }}.  In such a case, $\Phi_1\big(Y_0(x)\big)$ is analytic at $\hat{y}=Y_0(x_2)$. Moreover
$\gamma_2\big(x_2,Y_0(x_2)\big)\neq 0$.
From Lemma \ref{4-lem-5}, we have
\beqlb\label{4-88}
\Phi_2(\lambda \bar{u}_1)&&=-\frac{\gamma_1\big(\lambda \bar{u}_1,Y_0(\lambda \bar{u}_1)\big)\Phi_1\big(Y_0(\lambda \bar{u}_1)\big)}{\gamma_2\big(\lambda \bar{u}_1,Y_0(\lambda \bar{u}_1)\big)}.
\eeqlb
On the other hand, similar to \eqref{4-83}, we have
\beqlb\label{1-3}
Y_0(x_2)=Y_1(x_2)=x_\gamma \bar{u}_2.
\eeqlb
 From  \eqref{4-56}, \eqref{4-88}, \eqref{1-3} and \eqref{4-76}, we  have
 \beqlb\label{4-77}
 \lim_{\lambda\to \beta}(\beta-\lambda)\Phi_{\bar{u}}(\lambda)&&=\lim_{\lambda\to \beta}\frac{(\beta-\lambda)2(\gamma_2(\lambda \bar{u})-u_2\gamma(\lambda))\Phi_2(\lambda \bar{u}_1)}{\lambda (\beta-\lambda)}+\frac{2\gamma_1(\beta \bar{u})\Phi_1(\beta \bar{u}_2)}{\beta <\bar{u},\Sigma \bar{u}>}\nonumber
 \\&&=-\frac{2(\bar{u}_1r_{11}+\bar{u}_2r_{21})\Phi_1(\beta \bar{u}_2)}{\beta}+\frac{2\gamma_1(\beta \bar{u})\Phi_1(\beta \bar{u}_2)}{\beta <\bar{u},\Sigma \bar{u}>}.
 \eeqlb
It follows from Theorem \ref{TLT} that
 \beqlb\label{4-78}
\P\{<\bar{u},Z>>\lambda\}\sim K\exp\{-\beta\lambda\}.
\eeqlb

 \noindent{\bf {Subcase 3-3 $x_{\gamma}>z_0=\frac{x_2}{\bar{u}_1}$ with $\tilde{x}\neq x^*$}}. In this case, $\gamma(\beta)\neq 0$. Hence,
 \beqnn
 Y_0(x_2)=Y_1(x_2)\neq \beta \bar{u}_2.
 \eeqnn
 Therefore,
 \beqlb\label{4-92}
 \gamma_2(\beta \bar{u})\neq 0.
 \eeqlb
 It follows from Lemma \ref{thm1},\eqref{4-56} and \eqref{4-92} that
 \beqnn
 \lim_{\lambda\to \beta}\sqrt{\beta-\lambda}\Phi_{\bar{u}}(\lambda)&&=\lim_{\lambda\to \beta}\sqrt{\beta-\lambda}\frac{(\gamma_2(\lambda \bar{u})-u_2\gamma(\lambda))\Phi_2(\lambda \bar{u}_1)}{\gamma (\lambda)}
 \\&&=\frac{(\gamma_2(\beta \bar{u})-u_2\gamma(\beta))A_2(x_2)}{\gamma(\beta)}.
 \eeqnn
Then by Theorem \ref{TLT}, we have
\beqnn
 \P\{<\bar{u},Z>>\lambda\}\sim K \lambda^{-\frac{1}{2}} \exp\{-\beta\lambda\}.
\eeqnn

 \noindent{\bf{Subcase 3-4 $x_{\gamma}>z_0=\frac{x_2}{\bar{u}_1}$ with $\tilde{x}= x^*$}}.  In this case, we still have \eqref{4-92}. Hence,
from Lemma \ref{thm1}, \eqref{4-56} and \eqref{4-92}
  \beqnn
 \lim_{\lambda\to \beta}(\beta-\lambda)\Phi_{\bar{u}}(\lambda)&&=\lim_{\lambda\to \beta}(\beta-\lambda)\frac{\big(\gamma_2(\lambda \bar{u})-u_2\gamma(\lambda)\big)\Phi_2(\lambda \bar{u}_1)}{\gamma(\lambda)}
 \\&&=\frac{(\gamma_2(\beta \bar{u})-u_2\gamma(\beta))A_1(x_2)}{\gamma(\beta)}.
 \eeqnn
Then by Theorem \ref{TLT}, we have
\beqnn
 \P\{<\bar{u},Z>>\lambda\}\sim K \exp\{-\beta\lambda\}.
\eeqnn

 \noindent{\bf{Case 4 $\frac{x_2}{\bar{u}_1}<\min\{z_0,x_{\gamma}\}$ }}.  In this case $\beta=\frac{x_2}{\bar{u}_1}$. Moreover, since $x_2<\min\{x^*,\tilde{x}\}$,  we have \eqref{thm1-3}. Hence, \eqref{4-54} holds.  Moreover, in such a case,  we have
 \beqnn
 Y_0(\beta \bar{u}_1)=Y_1(\beta \bar{u}_1)\neq \beta \bar{u}_2.
 \eeqnn
 Therefore,
 \beqnn
 \gamma_2(\beta \bar{u})\neq 0.
 \eeqnn
Moreover, since $x_2<\min\{x^*,\tilde{x}\}$,  we can get that $\Phi_1(Y_0(\beta \bar{u}_1))$ is finite.  On the other hand, since $\beta_1>\beta_2$
\beqnn
\gamma_1(\beta \bar{u})\neq 0.
\eeqnn
 From \eqref{4-56} and \eqref{4-54}, we can get that
 \beqlb\label{4-79}
\lim_{\lambda\to \beta} \frac{\Phi_{\bar{u}}(\beta)-\Phi_{\bar{u}}(\lambda)}{K\sqrt{\beta-\lambda}}= 1.
 \eeqlb
Similar to \eqref{4-84}, by \eqref{4-79}, we get
\beqnn
\P\{<\bar{u},Z>>\lambda\}\sim K\lambda^{-\frac{3}{2}}\exp\{-\beta\lambda\}.
\eeqnn
\qed

\end{appendix}

\end{document}